\magnification=1200
\baselineskip=14pt

\font\tenbb=msbm10
\font\sevenbb=msbm7
\font\fivebb=msbm5
\newfam\bbfam
\textfont\bbfam=\tenbb \scriptfont\bbfam=\sevenbb
\scriptscriptfont\bbfam=\fivebb
\def\bb{\fam\bbfam}

\def\Ab{{\bb A}}
\def\Pb{{\bb P}}
\def\Qb{{\bb Q}}
\def\Zb{{\bb Z}}

\def\Ec{{\cal E}}
\def\Mc{{\cal M}}
\def\Oc{{\cal O}}
\def\Tc{{\cal T}}
\def\Wc{{\cal W}}
\def\Xc{{\cal X}}
\def\Yc{{\cal Y}}
\def\Zc{{\cal Z}}

\def\t{\theta}
\def\Lb{\Lambda}
\def\om{\omega}
\def\d{\delta}
\def\a{\alpha}
\def\vp{\varphi}
\def\s{\sigma}
\def\ve{\varepsilon}

\def\wh{\widehat}
\def\wt{\widetilde}
\def\part{\partial}
\def\ov{\overline}
\def\ra{\rightarrow}
\def\longra{\longrightarrow}

\def\ud{\underline}
\def\sbs{\subset}
\def\ot{\otimes}
\def\op{\oplus}
\def\ts{\times}
\def\bu{\bullet}
\def\ify{\infty}

\def\build#1_#2^#3{\mathrel{
\mathop{\kern 0pt#1}\limits_{#2}^{#3}}}

\def\hfl#1#2{\smash{\mathop{\hbox to 6mm{\rightarrowfill}}
\limits^{\scriptstyle#1}_{\scriptstyle#2}}}

\def\vfl#1#2{\llap{$\scriptstyle #1$}\left\downarrow
\vbox to 3mm{}\right.\rlap{$\scriptstyle#2$}}

\def\diagram#1{\def\normalbaselines{\baselineskip=0pt
\lineskip=10pt\lineskiplimit=1pt} \matrix{#1}}

\def\limind{\mathop{\oalign{lim\cr
\hidewidth$\longrightarrow$\hidewidth\cr}}}

\def\limproj{\mathop{\oalign{lim\cr
\hidewidth$\longleftarrow$\hidewidth\cr}}}


\vglue 1cm

\centerline{\bf Direct images in non-archimedean Arakelov 
theory}

\medskip

\centerline{\bf H. Gillet and C. Soul\'e}
\footnote{}{ The first author is supported by NSF  grant DMS-9801219}
\vglue 1cm

In this paper we develop a formalism of direct images for metrized 
vector bundles  in the context of 
the non-archimedean Arakelov theory introduced in our joint 
work 
[BGS] with S. Bloch, and we prove a Riemann-Roch-Grothendieck
theorem for this direct image. The new ingredient
 in the construction of the direct image is a non archimedean
 "analytic torsion current".

\smallskip

Let $K$ be the fraction field of a discrete valuation ring 
$\Lb$, and $X$ a smooth projective variety over $K$. In 
[BGS] 
we defined the codimension $p$ arithmetic Chow group of $X$ 
as the inductive limit
$$
\wh{\rm CH}^p (X) = \limind \, {\rm CH}^p (\Xc)
$$
of the Chow groups of the models $\Xc$
of $X$ over $\Lb$. Assuming 
resolution of singularities (cf.~1.1 below) we proved that 
these groups can also be defined as rational equivalence 
classes of pairs $(Z,g)$, where $Z$ is a codimension $p$ 
cycle 
on $X$, and $g$ is a ``Green current'' for $Z$. Here a 
``current'' is a projective system of cycle classes on the 
special fibers of all possible models of $X$. 
We have shown in 
[BGS] that many concepts and results 
in complex geometry and arithmetic 
intersection theory [GS1] have analogs in this context: 
differential forms, $\part \ov\part$-lemma, 
Poincar\'e-Lelong 
formula, intersection product, 
inverse and direct image maps in arithmetic Chow 
groups etc.

\smallskip

On the other hand, we defined a metrized vector bundle on 
$X$ 
to be a bundle $E$ on $X$, together with a bundle $E_{\Xc}$ 
on some model $\Xc$ of $X$ which restricts to $E$ on $X$. 
The 
theory of characteristic classes (resp. Bott-Chern secondary 
characteristic classes) for hermitian vector bundles on 
arithmetic varieties [GS2] is replaced here by 
characteristic 
classes with values in the Chow groups of
$\Xc$ (resp. the Chow groups of $\Xc$ 
with supports in its special fiber) ([BGS], (1.9), and \S~2 
below). These classes
are contravariant for maps of varieties over 
$K$.

\smallskip

However, we were not able in [BGS] to define direct images 
of 
metrized vector bundles. Recall that in Arakelov geometry, 
if 
$f : X \ra Y$ is a map of varieties over $\Zb$ which is 
smooth on the set of complex points of $X$, and if $\ov E$ 
is 
an hermitian vector bundle on $X$, once we choose a metric on 
$Tf$, the $L^2$-metric on the determinant
line bundle $det(Rf_* \, E)$ need {\it 
not} be smooth in general. For this reason, following an 
idea 
of Quillen [Q], one is led to modify the $L^2$-metric 
on the determinant line bundle by 
multiplying it by the Ray-Singer analytic torsion
of the Dolbeault complex, which results in a smooth metric.
One of
the key features of the Quillen metric, is that it gives a Riemann-Roch
formula for the first Chern class of the determinant line bundle which is
an equality of forms.
More generally,  if one chooses a complex of vector
bundles $F.$, a quasi-isomorphism 
       $ F. \ra Rf_*E$  
and hermitian metrics on all the $F_i$'s, one can define 
a form $\theta$ on the complex points of $Y$,
called the higher analytic torsion 
([GS3],  [BK]), which is well defined
up to boundaries, and such that
$ - dd^c(\theta)$ is equal to the difference between the 
Chern character form of $F.$ and the direct image
of the product of  the Chern character
  form of $\ov E$
  with the Todd form of $Tf$. This form $\theta$ is the
 key ingredient
 when defining direct images for the "arithmetic Grothendieck groups"
[GS2] [GS3].

\smallskip

In the non-archimedean case we face a similar difficulty. 
Assume $f : X \ra Y$ is a morphism of projective varieties 
over $K$, induced by a map of models $\ud f : \Xc \ra \Yc$. 
Let $E_{\Xc}$ be a vector bundle on $\Xc$, with restriction 
$E$ to $X$. A natural candidate for a (non-archimeadean)
metric 
on $Rf_* \, E$ is then the complex of vector bundles
 $R\ud{f}_* \, E_{\Xc}$ on 
$\Yc$. But this choice
will {\it not}, in general,
 be compatible with changes of 
models for both $X$ and $Y$. Indeed, consider a commutative 
diagram of models
$$
\diagram{
\Xc' &\hfl{\pi}{} &\Xc \cr
\vfl{\ud{f}'}{} &&\vfl{}{\ud{f}} \cr
\Yc' &\hfl{}{p} &\Yc \cr
}
$$
where both $\ud f$ and $\ud f'$ induce $f : X \ra Y$. The 
canonical map
$$
Lp^* \, R \ud{f}_* \, E_{\Xc} \ra R \ud{f}'_* \, L \pi^* \, 
E_{\Xc}
$$
need not be an isomorphism. We are led to use the Chern 
character with supports of the cone of this map to define 
the 
$\Yc'$-component $\t_{\Yc'}$ of the ``higher analytic 
torsion'' $\t$, which is a current on $Y$ (see (50) and 
Prop.~4 
for a precise definition).

\smallskip

We then define a Grothendieck group ${\build {K_0}_{}^{\vee}} 
(X)$, generated by triples $(E,h,\eta)$ where $E$ is a 
bundle 
on $X$, $h$ is a metric on $E$ and $\eta$ is a sum of
currents of all degrees on $X$. The relations in 
${\build {K_0}_{}^{\vee}} 
(X)$ come from exact sequences of vector bundles on $X$ 
(\S~2.6, (37)). By imposing that $\eta$ be smooth 
(i.e. $\eta$ consists of an inductive system of cycle
classes and not only a projective one,
see [BGS] and 1.2. below),
we also 
define a subgroup $\wh{K}_0 (X) \sbs {\build {K_0}_{}^{\vee}} (X)$. Now 
let $f : X \ra Y$ be any morphism, and choose a metric on 
the 
relative tangent complex $Tf$. We attach to these data a 
direct image morphism
$$
f_* : {\build {K_0}_{}^{\vee}} (X) \ra {\build {K_0}_{}^{\vee}} (Y) \, .
$$
If $\ud f : \Xc \ra \Yc$ is a map of models inducing $f$, if 
$E$ is a bundle on $X$, if the metric $h$ on $E$ (resp. the 
metric on $Tf$) is defined by a bundle $E_{\Xc}$ on $\Xc$ 
(resp. by $T \ud f$), and if $R^q \ud{f}_* \, E_{\Xc} = 0$, 
$q > 0$, the direct image $f_* (E,h,0)$ is the class of 
$(f_* 
\, E ,  \ud{f}_* \, E_{\Xc} , \t)$, where $\t$ is the higher 
analytic torsion of $E_{\Xc}$ (Prop.~4, Th.~1).

\smallskip

When $f$ is flat, $f_*$ maps $\wh{K}_0 (X)$ into $\wh{K}_0 
(Y)$ and a Riemann-Roch-Grothendieck theorem holds for Chern 
characters with values in $\wh{\rm CH}^* \ot \Qb$ (Th.~1, i) 
and Th.~2, ii)). This is not so surprising, as it follows 
from the definition of $\t$ and the 
Riemann-Roch-Grothendieck theorem with values in Chow groups of projective 
schemes over $\Lb$. What is more involved is, first, to show that the family $\t 
= (\t_{\Yc'})$ does define a current on $Y$ (a form when $f$ is flat) (Prop.~4) and, second, to check that the expected anomaly formulae for the 
change of metrics on either $E$ or $Tf$ are true in our case 
(Th.~1, (64) and (66)).
These facts rely upon the vanishing of the direct image of the relative 
Todd class 
with support  of  birational maps (Prop.3 ii)). This key lemma 
is itself a consequence  of the proof by Franke of a refined Riemann-Roch 
formula conjectured by Saito 
(this proof of Franke [Fr] remains unfortunately unpublished);
we also give an alternative proof for birational maps
for which the strong factorization conjecture holds.
Finally we prove that our direct
image is such that $(gf)_* = g_* f_*$ and 
that a projection formula is valid 
(Th.~2, i) and iii)).

\smallskip

The paper is organized as follows. In \S~1 we review the 
arithmetic intersection theory of [BGS], and we also 
introduce a group ${\build {{\rm CH}^p}_{}^{\vee}} (X)$ containing $\wh{\rm 
CH}^p (X)$ which is always covariant (this definition was 
inspired by similar definitions by Burgos [B] and Zha [Z]). 
In \S~2 we develop the theory of secondary characteristic 
classes for metrized complexes of vector bundles and we 
define arithmetic Grothendieck groups $\wh{K}_0 (X) \sbs 
{\build {K_0}_{}^{\vee}} (X)$ together with Chern characters from 
${\build {K_0}_{}^{\vee}} (X)$ (resp. $\wh{K}_0 (X)$) to $\build \op_{p 
\geq 0}^{} {\build {{\rm CH}^p}_{}^{\vee}} (X) \ot \Qb$ (resp. $\build 
\op_{p \geq 0}^{} \wh{\rm CH}^p (X) \ot \Qb$). In \S~3 we 
prove several facts about the secondary Todd classes of 
birational maps of models. In \S~4, after defining the 
higher 
analytic torsion currents (Prop.~4), we define the direct 
images $f_*$ on ${\build {K_0}_{}^{\vee}}$, and we give properties of 
this 
map, including a Riemann-Roch-Grothendieck theorem with 
values in ${\build {\rm CH}_{}^{\vee}}$ (Th.~1 and Th.~2).

\smallskip

We thank S.~Bloch
for useful discussions. We are also 
grateful to the Newton Institute where most of this work was 
done.

\bigskip

\noindent {\bf 1. Cycles}

\smallskip

\noindent {\bf 1.1.} We first recall some definitions and 
results in ``non-archimedean Arakelov theory'', from our 
joint work with S.~Bloch [BGS], to which we refer the reader 
for more details.

\smallskip

Let $\Lb$ be an excellent discrete valuation ring with 
quotient field $K$ and residue field $k$. Let $X$ be a 
smooth 
projective scheme over ${\rm Spec} (K)$. Let $\Xc$ be a 
regular scheme, projective and flat over ${\rm Spec}(\Lb)$, 
together with an isomorphism of its generic fiber
$\Xc_K$ 
with $X$. We denote by $\Xc_0 = \Xc \build \ts_{{\rm Spec} 
(\Lb)}^{} {\rm Spec} (k)$ the special fiber of $\Xc$
and by $i : 
\Xc_0 \ra \Xc$ its closed immersion into $\Xc$. We shall say 
that $\Xc$ is a {\it model} of $X$ when $\Xc_0^{\rm red}$ is 
a divisor with normal crossings.

\smallskip

A {\it map of varieties} $X \ra Y$ will mean a morphism over 
$K$ between smooth projective schemes over $K$. When $\Xc$ 
and $\Yc$ are models of $X$ and $Y$ respectively, a {\it map 
of models} $\Xc \ra \Yc$ is any morphism defined over $\Lb$. 
A map of models is {\it good} when it is the composite of 
blow ups with integral regular centers meeting the reduced 
special fiber normally.

\smallskip

We shall assume that axioms (M1) and (M2) of [BGS] (1.1) 
hold. Axiom (M1) says that, given any scheme
$\Xc$, projective 
and flat over $\Lb$, with smooth generic fiber $X$, there 
exists a model $\Xc'$ of $X$ and a morphism $\Xc' \ra \Xc$ 
over $\Lb$. Axiom (M2) says that, given two models $\Xc$ 
and $\Xc'$ of $X$, there exists a third model $\Xc''$, a map 
of models $\Xc \ra \Xc'$, and a good map of models $\Xc \ra 
\Xc''$.

By Hironaka [H], these axioms are satisfied when $\Lb$ is a 
localization of an algebra of finite type over a field of 
characteristic zero.

\smallskip

We write $\Mc (X)$ for the category of models of $X$.

\medskip

\noindent {\bf 1.2.} Under the assumption of 1.1, let $\Xc$ 
be a model of $X$ and let $p \geq 0$ be an integer. Denote by 
${\rm CH}_p (\Xc_0)$ the Chow homology group of dimension 
$p$ algebraic cycles on $\Xc_0$ modulo rational equivalence, 
and by ${\rm CH}^p (\Xc_0)$ the Chow cohomology group of 
codimension $p$ of $\Xc_0$, i.e. the bivariant group ${\rm 
CH}^p (\Xc_0 \build \longra_{}^{\rm id} \Xc_0$) of [Fu], 
17.1.

\smallskip

When $X$ is fixed, any morphism $\pi : \Xc' \ra \Xc$ between 
models of $X$ induces both direct images
$$
\displaylines{
\pi_* : {\rm CH}_p (\Xc'_0) \ra {\rm CH}_p (\Xc_0) \cr
\pi_* : {\rm CH}^p (\Xc'_0) \ra {\rm CH}^p (\Xc_0) \cr
}
$$
and inverse images
$$
\displaylines{
\pi^* : {\rm CH}_p (\Xc_0) \ra {\rm CH}_p (\Xc'_0) \cr
\pi^* : {\rm CH}^p (\Xc_0) \ra {\rm CH}^p (\Xc'_0) \cr
}
$$
([BGS] (1.4)). The projection formula implies
$$
\pi_* \, \pi^* = {\rm id} \, . \leqno (1)
$$
If $d$ is the dimension of $\Xc_0$ over $k$, we may consider 
the inductive limits with respect to $\pi^*$:
$$
\eqalign{
A_{\rm closed}^{pp} (X) &:= \limind_{\Mc (X)} {\rm CH}^p 
(\Xc_0) \cr
\wt{A}^{pp} (X) &:= \limind_{\Mc (X)} {\rm CH}_{d-p} (\Xc_0) 
\, , \cr
}
$$
as well as the projective limits under $\pi_*$:
$$
\eqalign{
D_{\rm closed}^{pp} (X) &:= \limproj_{\Mc (X)} {\rm CH}^p 
(\Xc_0) \cr
\wt{D}^{pp} (X) &:= \limproj_{\Mc (X)} {\rm CH}_{d-p} 
(\Xc_0) \, . \cr
}
$$
By analogy with Arakelov theory [GS1] these groups are 
called, respectively, closed $(p,p)$-forms, $(p,p)$-forms 
modulo boundaries, closed $(p,p)$-currents, and 
$(p,p)$-currents modulo boundaries.

From (1) it follows that there are 
canonical inclusions
$$
A_{\rm closed}^{pp} (X) \sbs D_{\rm closed}^{pp} (X)
$$
and
$$
\wt{A}^{pp} (X) \sbs \wt{D}^{pp} (X)
$$
of forms into currents.

\smallskip

We shall also denote by $A_{\rm closed}^{pp} (X)_{\Qb}$, 
$\wt{A}^{pp} (X)_{\Qb}$,~$\ldots$ the tensor products 
$A_{\rm closed}^{pp} (X) \build \ot_{\Zb}^{} \Qb$, 
$\wt{A}^{pp} (X) \build \ot_{\Zb}^{} \Qb$,~$\ldots$ 
Furthermore we let
$$
A_{\rm closed} (X)_{\Qb} = \build \op_{p \geq 0}^{} A_{\rm 
closed}^{pp} (X)_{\Qb} \, ,
$$
and we define similarly $\wt{A} (X)_{\Qb}$ etc~$\ldots$

\medskip

\noindent {\bf 1.3.} Given a model $\Xc$, there is a 
morphism
$$
i^* \, i_* : {\rm CH}_{d-p} (\Xc_0) \ra {\rm CH}^{p+1} 
(\Xc_0) \, ,
$$
obtained by composing the direct image in Chow homology
$$
i_* : {\rm CH}_{d-p} (\Xc_0) \ra {\rm CH}_{d-p} (\Xc) \, ,
$$
the Poincar\'e duality
$$
{\rm CH}_{d-p} (\Xc) \simeq {\rm CH}^{p+1} (\Xc) \, ,
$$
and the restriction map in Chow cohomology
$$
i^* : {\rm CH}^{p+1} (\Xc) \ra {\rm CH}^{p+1} (\Xc_0) \, .
$$
When $\Xc$ varies, these maps are compatible with $\pi_*$ 
and they induce a morphism
$$
dd^c : \wt{D}^{pp} (X) \ra D^{p+1 , p+1} (X)
$$
on projective limits. Using resolution of singularities, one 
gets the following result ([BGS], Th.~2.3.1):

\medskip

\noindent {\bf Proposition 1.}

\item{i)} {\it A current $g \in \wt{D}^{pp} (X)$ lies in 
$\wt{A}^{pp} (X)$ if and only if $dd^c (g)$ lies in the 
subgroup $A_{\rm closed}^{p+1 , p+1} (X)$ of $D_{\rm 
closed}^{p+1 , p+1} (X)$.}

\item{ii)} {\it The kernel (resp. cokernel) of $dd^c$ 
coincides with the kernel (resp. cokernel) of $i^* \, i_*$ 
on any model of $X$.}

\medskip

\noindent {\bf 1.4.} Let $Y \sbs X$ be a closed integral 
subvariety of codimension $p$. For any model $\Xc$ we may 
consider the Zariski closure $\ov Y$ of $Y$ in $\Xc$, and 
the restriction $i^* [\ov Y] \in {\rm CH}^p (\Xc_0)$ of its 
cycle class on $\Xc$. These are compatible with $\pi_*$ and 
we get this way a closed current
$$
\d_Y = (i^* [\ov Y]) \in D_{\rm closed}^{pp} (X) \, . 
$$
When $Z = \ \build \sum_{\a}^{} n_{\a} \, Y_{\a} \in Z^p 
(X)$ is any algebraic cycle of codimension $p$ on $X$, we 
let
$$
\d_Z = \sum_{\a} n_{\a} \, \d_{Y_{\a}} \, .
$$
A {\it Green current} for $Z$ is any current $g \in 
\wt{D}^{p-1,p-1} (X)$ such that $dd^c (g) + \d_Z$ lies in 
$A_{\rm closed}^{pp} (X)$.

\smallskip

For example, let $W \sbs X$ be a closed integral subvariety 
of codimension $p-1$ and $f \in k(W)^*$ a non trivial 
rational function on $W$. Let ${\rm div} (f)$ be the divisor 
of $f$ on $W$, viewed as a codimension $p$ cycle on $X$, let 
$\ov{{\rm div} (f)}$ be its Zariski closure on a model 
$\Xc$, and let ${\rm div}_{\Xc} (f)$ be the divisor of $f$ 
on the Zariski closure of $W$ in $\Xc$. Consider the 
difference
$$
{\rm div}_{\nu} (f)_{\Xc} = \ov{{\rm div} (f)} - {\rm 
div}_{\Xc} (f) \, .
$$
The family $-{\rm div}_{\nu} (f) = (-{\rm div}_{\nu} 
(f)_{\Xc})$ is then a Green current for the cycle ${\rm div} 
(f)$ ([BGS], (3.1)).

\medskip

\noindent {\bf 1.5.} For any $p \geq 0$, the {\it arithmetic 
Chow group} $\wh{\rm CH}^p (X)$ is defined as the quotient 
of the abelian group of pairs $(Z,g)$, where $Z \in Z^p (X)$ 
is a codimension $p$ algebraic cycle on $X$ and $g$ is a 
Green current for $X$, by the subgroup generated by 
the set of all 
pairs $({\rm div} (f), -{\rm div}_{\nu} (f))$, where $f$ is 
a non zero rational function on a codimension $(p-1)$ closed 
integral subvariety in $X$.

\smallskip

Let us also introduce the group ${\build {{\rm CH}^p}_{}^{\vee}} (X)$, 
equal to the quotient of
 $Z^p (X) \op \wt{D}^{p-1,p-1} (X)$ by the group 
generated by pairs $({\rm div} (f) , -{\rm div}_{\nu} (f))$ 
as above (compare [B] and [Z]). Clearly there is an inclusion
$$
\wh{\rm CH}^p (X) \sbs {\build {{\rm CH}^p}_{}^{\vee}} (X) \, .
$$
Let
$$
\om : {\build {{\rm CH}^p}_{}^{\vee}} (X) \ra D_{\rm closed}^{pp} (X)
$$
be the morphism sending the class of $(Z,g)$ to $dd^c (g) + 
\d_Z$ (this kills the relations in ${\build {{\rm CH}^p}_{}^{\vee}} (X)$, 
cf. [BGS] Prop.~(3.1.1)). 
The subgroup $\wh{\rm CH}^p (X)$ 
consists of those $x$ in ${\build {{\rm CH}^p}_{}^{\vee}} (X)$ such that 
$\om (x)$ lies in the subgroup
$A_{\rm closed}^{pp} (X)$ of $D_{\rm closed}^{pp} (X)$.

\smallskip

We also denote by
$$
a : \wt{D}^{p-1,p-1} (X) \ra {\build {{\rm CH}^p}_{}^{\vee}} (X)
$$
the morphism sending $\eta$ to the class of $(0,\eta)$. 
By Proposition 1 i), $\eta$ lies in $\wt{A}^{p-1,p-1} (X)$ 
if and only if 
its image $a(\eta)$ lies in $\wh{\rm CH}^p (X)$.

\smallskip

As shown in [BGS] Th.~3.3.3, there is a canonical 
isomorphism
$$
\limind_{\Mc (X)} {\rm CH}^p (\Xc) \build \longra_{}^{\sim} 
\wh{\rm CH}^p (X) \leqno (2)
$$
and (taking inverse limit in the diagram of op.cit., p.~461) 
it extends to an isomorphism
$$
\limproj_{\Mc (X)} {\rm CH}^p (\Xc) \build \longra_{}^{\sim} 
{\build {{\rm CH}^p}_{}^{\vee}} (X) \, . \leqno (3)
$$
Given $\eta = (\eta_{\Xc}) \in \wt{D}^{p-1,p-1} (X)$ this 
isomorphism sends $(i_* \, \eta_{\Xc})$ to $a(\eta)$, and if 
$Z \in Z^p (X)$ it sends the projective system of Zariski 
closures of $Z$ in $\Xc$ to the class of $(Z,0)$.

\medskip

\noindent {\bf 1.6.} Let $f : X \ra Y$ be a map of 
varieties. We know from
[BGS] 1.6 that 
forms are 
contravariant and currents are covariant.
Furthermore, it follows from (2) that $f$ induces pull-back 
morphisms
$$
f^* : \wh{CH}^p (Y) \ra \wh{CH}^p (X) 
$$
and from (3) we get direct image morphisms
$$
f_* : {\build {{\rm CH}^p}_{}^{\vee}} (X) \ra {\build {\rm CH}_{}^{\vee}}^{p-\d} (X) \, 
,
$$
where $\d$ is the relative dimension of $X$ over $Y$. We may 
also describe $f_*$ as mapping the class of $(Z,g)$ to the 
class of $(f_* (Z), f_* (g))$, where $f_* (Z)$ is the usual 
direct image of the cycle $Z$ [Fu]. Given 
two maps of varieties $f : X \ra Y$ and $h : Y \ra Z$, we 
have $(hf)_* = h_* \, f_*$ and $(hf)^* = f^* \, h^*$.

\smallskip

Assume furthermore that $f$ is flat. Then, as was shown in 
[BGS], Thms.~(4.1.1) and (4.2.1), the morphism $f_*$ maps 
forms to forms and respects $\wh{\rm CH}$:
$$
\eqalign{
f_* (\wt{A}^{pp} (X)) &\sbs \wt{A}^{p-\d , p-\d} (X) \, , 
\cr
f_* (A_{\rm closed}^{pp} (X)) &\sbs A_{\rm closed}^{p-\d , 
p-\d} (X) \, , \cr
f_* (\wh{\rm CH}^p (X)) &\sbs \wh{\rm CH}^{p-\d} (X) \, .
}
$$

\smallskip

\noindent {\bf 1.7.} From formula (2) we also deduce a graded 
ring structure
$$
\wh{\rm CH}^p (X) \ot \wh{\rm CH}^q (X) \ra \wh{\rm 
CH}^{p+q} (X) \leqno (4)
$$
on arithmetic Chow groups. From (3) and the projection 
formula $\pi_* (x \, \pi^* (y)) = \pi_* (x) y$ for any map 
$\pi : \Xc' \ra \Xc$ between models of $X$, $x \in {\rm 
CH}^p (\Xc')$, $y \in {\rm CH}^q (\Xc)$, we deduce a pairing
$$
{\build {{\rm CH}^p}_{}^{\vee}} (X) \ot \wh{\rm CH}^q (X) \ra 
{\build {{\rm CH}^{p+q}}_{}^{\vee}} (X) \leqno (5)
$$
extending (4), and turning ${\build {\rm CH}_{}^{\vee}} (X)$ into a graded 
module on $\wh{\rm CH} (X)$.

\smallskip

When $f : X \ra Y$ is a map of varieties, the formulae
$$
f^* (xy) = f^* (x) \, f^* (y)
$$
and
$$
f_* (x \, f^* (y)) = f_* (x) \, y
$$
hold when $x$ lies in $\wh{\rm CH}^p (Y)$ (resp. $\check{\rm 
CH}^p (X)$) and $y$ lies in $\wh{\rm CH}^q (Y)$.

\smallskip

Similar facts are true for pairings between forms and 
currents ([BGS] 1.5 and Proposition (1.6.2)).

\bigskip

\noindent {\bf 2. Vector bundles}

\smallskip

\noindent {\bf 2.1.} We keep the notation of sections 1.1 or 
1.2. Let $E$ be a vector bundle on $X$. A {\it metric} on 
$E$ is determined by the choice of a vector bundle $h = 
E_{\Xc}$ on some model $\Xc$ of $X$, together with an 
isomorphism $E \simeq E_{\Xc \mid X}$ of $E$ with the 
restriction of $E_{\Xc}$ to $X$. By convention,
given any map $\pi : \Xc' 
\ra \Xc$ between models of $X$, $\pi^* \, E_{\Xc}$ defines 
the same metric as $E_{\Xc}$ (see [BGS] (1.9.1)). Notice 
that (by [RG] Part I, Th.~5.2.2 together with [Ma] Th.~7.10) 
any bundle $E$ on $X$ has a metric. Furthermore, given any 
morphism of bundles $u : E \ra F$ and arbitrary metrics on 
$E$ and $F$, there exists a model $\Xc$ and a (unique) 
morphism of bundles $E_{\Xc} \ra F_{\Xc}$ which induces $u$ 
on $X$ and such that $E_{\Xc}$ (resp. $F_{\Xc}$) induces the 
given metric on $E$ (resp. $F$) (by uniqueness this is a local 
problem on $\Xc$, therefore we can assume that both 
$E_{\Xc}$ and $F_{\Xc}$ are trivial; after blowing up $\Xc$, 
all the coefficients of the matrix $u$ of rational functions 
extend to some model of $X$).

\smallskip

Let $D^b (\Xc)$ (resp. $D^b (X)$) be the derived category of 
bounded complexes of vector bundles on $\Xc$ (resp. $X$). 
Since $\Xc$ (resp. $X$) is noetherian and regular, this is 
equivalent to the derived category of bounded complexes of 
coherent sheaves on $\Xc$ (resp. $X$). We also denote 
by $D_{\Xc_0}^b (\Xc)$ the full subcategory 
of $D^b (\Xc)$ consisting of complexes which are acyclic 
outside $\Xc_0$.

\smallskip

Given $K \in D^b 
(X)$, a {\it metric} on $K$ is determined by the choice of 
$K_{\Xc} \in D^b (\Xc)$ together with an isomorphism $K 
\simeq K_{\Xc \mid X}$. Again, given $\pi : \Xc' \ra \Xc$, 
$\pi^* \, K_{\Xc}$ defines the same metric as $K_{\Xc}$.

\medskip

\noindent {\bf 2.2.} Let $\phi = {\rm ch}$ or ${\rm Td}$ be 
the Chern character or the Todd genus. Given any metrized 
vector bundle $\ov E = (E,h)$ on $X$, we can attach to $\ov 
E$ a closed form
$$
\phi (\ov E) \in A_{\rm closed} (X)_{\Qb} \, .
$$
If $E_{\Xc}$ is an extension of $E$ defining $h$, this form 
$\phi (\ov E)$ is defined as the image in the direct limit 
of Chow cohomology groups of
$$
\phi (i^* \, E_{\Xc}) \in {\rm CH}^{\bu} (\Xc_0)_{\Qb} = \ 
\build \op_{p \geq 0}^{} {\rm CH}^p (\Xc_0)_{\Qb} \, .
$$
The Chern character is such that
$$
{\rm ch} (\ov E \op \ov F) = {\rm ch} (\ov E) + {\rm ch} 
(\ov F)
$$
and
$$
{\rm ch} (\ov E \ot \ov F) = {\rm ch} (\ov E) \, {\rm ch} 
(\ov F) \, ,
$$
while the Todd class is multiplicative
$$
{\rm Td} (\ov E \op \ov F) = {\rm Td} (\ov E) \, {\rm Td} 
(\ov F) \, .
$$
Here, given $\ov E = (E,E_{\Xc})$ and $\ov F = 
(F,F_{\Xc'})$, their sum $\ov E \, \op \, \ov F$ is defined as $E \, \op \, F$ 
extended to $\pi^* \, E_{\Xc} \, \op \, (\pi')^* \, F_{\Xc'}$ on any model 
$\Xc''$ with maps $\pi : \Xc'' \ra \Xc$ and $\pi' : \Xc'' \ra \Xc'$.

\smallskip

From (2) we may also define a class
$$
\wh{\phi} (\ov E) \in \wh{\rm CH}^{\bu} (X)_{\Qb}
$$
such that $\om (\wh{\phi} (\ov E)) = \phi (\ov E)$. This is 
just the image of $\phi (E_{\Xc}) \in {\rm CH}^{\bu} 
(\Xc)_{\Qb}$ in the inductive limit.

\smallskip

When $K \in D^b (X)$ is equipped with a metric we may also 
define $\phi (\ov K) \in A_{\rm closed} (X)_{\Qb}$ and 
$\wh{\phi} (\ov K) \in \wh{\rm CH} (X)_{\Qb}$ such that $\om 
(\wh{\phi} (\ov K)) = \phi (\ov K)$. If $C^{\bu}$ is a 
bounded complex of bundles on $\Xc$ representing the metric 
on $K$, we just let $\wh{\rm ch} (\ov K)$ be the image of
$$
{\rm ch} (C^{\bu}) := \sum_{i \geq 0} (-1)^i \, {\rm ch} 
(C^i)
$$
and $\wh{\rm Td} (\ov K)$ be the image of
$$
{\rm Td} (C^{\bu}) := \prod_{i \geq 0} {\rm Td} 
(C^i)^{(-1)^i}
$$
in the inductive limit 
of Chow groups of models of $X$. This is 
well defined since ${\rm ch}$ (resp. ${\rm Td}$) is additive 
(resp. multiplicative) on exact sequences of vector bundles 
on $\Xc$.

\medskip

\noindent {\bf 2.3.} Let $X$ be as above and let $\Xc$ be a 
model of $X$. Consider a bounded acyclic complex
$$
(E^{\bu} , d) = (E^0 \build \longra_{}^{d} E^1 \build 
\longra_{}^{d} \cdots \build \longra_{}^{d} E^k)
$$
of bundles over $X$. For every $n \geq 0$, let $E_{\Xc}^n$ 
be a vector bundle on $\Xc$ restricting to $E^n$ on $X$. Let 
$\phi = {\rm ch}$ or ${\rm Td} - 1$ and denote by $\phi 
(E_{\Xc}^{\bu})$ the elements
$$
{\rm ch} (E_{\Xc}^{\bu}) = \sum_{n=0}^k (-1)^n \, {\rm ch} 
(E_{\Xc}^n) \leqno (6)
$$
and
$$
({\rm Td} - 1) (E_{\Xc}^{\bu}) = (\prod_{n=0}^k {\rm Td} 
(E_{\Xc}^n)^{(-1)^n} )- 1 \leqno (7)
$$
in ${\rm CH} (\Xc)_{\Qb}$.

\medskip

\noindent {\bf Proposition 2.} {\it Under these assumptions 
there exists a unique class
$$
\phi_{\Xc_0} (E_{\Xc}^{\bu}) \in {\rm CH}_{\bu} 
(\Xc_0)_{\Qb}
$$
with the following three properties:}

\item{i)}{\it One has}
$$i_* \, \phi_{\Xc_0} (E_{\Xc}^{\bu}) = \phi 
(E_{\Xc}^{\bu}).$$

\item{ii)} {\it Let $f : \Yc \ra \Xc$ be any map of models. 
Then}
$$
f^* \, \phi_{\Xc_0} (E_{\Xc}^{\bu}) = \phi_{\Xc_0} (f^* 
(E_{\Xc}^{\bu})) \, .
$$

\item{iii)} {\it Assume that, for all $n \geq 0$, the 
differential $d : E^n \ra E^{n+1}$ extends to $d_{\Xc} : 
E_{\Xc}^n \ra E_{\Xc}^{n+1}$ and that $(E_{\Xc}^{\bu} , 
d_{\Xc})$ is acyclic on $\Xc$. Then}
$$
\phi_{\Xc_0} (E_{\Xc}^{\bu}) = 0 \, .
$$
\smallskip

{\it Furthermore}

\item{iv)} {\it Let
$$
0 \ra S^{\bu} \ra E^{\bu} \ra Q^{\bu} \ra 0
$$
be an exact sequence of bounded acyclic complexes on $X$. 
Assume $S_{\Xc}^n$, $E_{\Xc}^n$, $Q_{\Xc}^n$ are bundles on 
$\Xc$ which restrict to $S^n$, $E^n$ and $Q^n$ 
respectively on $X$, $n \geq 0$. 
Then the following identities hold:}
$$
{\rm ch}_{\Xc_0} (E_{\Xc}^{\bu}) = {\rm ch}_{\Xc_0} 
(S_{\Xc}^{\bu}) + {\rm ch}_{\Xc_0} (Q_{\Xc}^{\bu}) \, , 
\leqno (8)
$$
$$
\leqalignno{
({\rm Td}-1)_{\Xc_0} (E_{\Xc}^{\bu}) = \ &({\rm Td}-1)_{\Xc_0} 
(S_{\Xc}^{\bu}) \cdot (i^* \, {\rm Td} (Q_{\Xc}^{\bu})) + 
({\rm Td}-1)_{\Xc_0} (Q_{\Xc}^{\bu}) &(9) \cr
= \ &({\rm Td}-1)_{\Xc_0} (S_{\Xc}) + i^* ({\rm Td} 
(S_{\Xc}^{\bu})) \cdot ({\rm Td}-1)_{\Xc_0} (Q_{\Xc}^{\bu}) 
\, . \cr
}
$$

\item{v)} {\it If $F_{\Xc}$ is any bundle on $\Xc$,}
$$
{\rm ch}_{\Xc_0} (E_{\Xc}^{\bu} \ot F_{\Xc}) = {\rm 
ch}_{\Xc_0} (E_{\Xc}^{\bu}) \, i^* \, {\rm ch} (F_{\Xc}) \, 
. \leqno (10)
$$

\item{vi)} {\it If $E_{\Xc}^{\bu} [1]$ denotes the shift by one
of $E_{\Xc}^{\bu}$, the  following two equalities hold:\it}
$$
{\rm ch}_{\Xc_0} (E_{\Xc}^{\bu} [1]) = - {\rm ch}_{\Xc_0} 
(E_{\Xc}^{\bu}) \leqno (11)
$$
{\it and}
$$
({\rm Td} - 1)_{\Xc_0} (E_{\Xc}^{\bu} [1]) = - ({\rm Td} - 
1)_{\Xc_0} (E_{\Xc}^{\bu}) \, {\rm Td} (E_{\Xc}^{\bu})^{-1} \, . \leqno 
(12)
$$

\smallskip

\noindent {\it Proof.} 
We use  the 
Grassmann-graph construction [BFM] (see also [Fu] \S 18.1 
, and [GS4] \S 1 for a variant of this construction).
 Let 
$e_n$ be the rank of $E_{\Xc}^n$, consider the Grassmann variety 
$G_n = {\rm Grass}_{e_n} (E_{\Xc}^n \op E_{\Xc}^{n+1})$, $n 
\geq 0$, and let
$$
G = G_0 \build \ts_{\Xc}^{} G_1 \build \ts_{\Xc}^{} \cdot 
\build \ts_{\Xc}^{} G_{k-1} \, .
$$
The acyclic complex $(E^{\bu} , d)$ defines a section
$$
\vp : X \build \ts_{K}^{} \Pb^1 
\ra G \build \ts_{\Lb}^{} \Pb^1
$$
of the projection
$$
G \build \ts_{\Lb}^{} \Pb^1 \ra \Xc \build \ts_{\Lb}^{} 
\Pb^1 \, 
$$
on the generic fiber $X \build \ts_{K}^{} \Pb^1$
(given by the graphs of the maps $\lambda d$
at the point $(x,\lambda) \in X \build \ts_{K}^{} \Ab^1$). 
We let $\Wc_1$ be the Zariski closure of $\vp (X \build 
\ts_{K}^{} \Pb^1)$ in $G \build \ts_{\Lb}^{} \Pb^1$, and 
$\Wc \ra \Wc_1$ a resolution of $\Wc_1$. The scheme
 $\Wc$ is a model of  $X \build \ts_{K}^{} \Pb^1$.
For each $n\geq 0$, the tautological bundle of rank 
$e_n$ on $G$ defines a vector bundle ${E}_{\Wc}^{n}$
on $\Wc$, and there exists an acyclic complex 
of vector bundles $(\wt{E}^{\bu} , \wt d)$ on $X \build \ts_{K}^{} \Pb^1$
where, for each $n\geq 0$, $\wt{E}^{n}$
is the restriction of ${E}_{\Wc}^{n}$ to $X \build \ts_{K}^{} \Pb^1$.

Let $\Wc_0$ (resp. $\Wc_{\ify}$) be the  
 Zariski closure $X \ts \{ 0 \}$ (resp. $X \ts \{ \ify \}$) in $\Wc$.
Denote by
$$
j_0 : \Xc =  \Wc_0 \ra \Wc
$$
the inclusion and by
$$
j_{\ify} : \Xc' \ra \Wc
$$
the composite of a resolution
of singularities $\Xc' \ra  \Wc_{\ify}$
with the inclusion
$ \Wc_{\ify} \ra \Wc$.
One has $j_0^* ({E}_{\Wc}^{\bu}) = E_{\Xc}^{\bu}$
and there exists a split acyclic complex
 $(j_{\ify}^* ({E}_{\Wc}^{\bu}), d_{\ify})$
on $\Xc'$ which restricts to
$j_{\ify}^*(\wt{E}^{\bu} , \wt d)$ on 
$X \ts \{ \ify \}$
([Fu], proof of Lemma 18.1). 

The standard parameter $z$ on $\Pb^1$ 
defines a rational function on $\Wc$, hence a class
$$
\ell (z) =  {\rm div}_{\nu} \, (z)_{\Wc}
$$
such that
$$
i_* \, \ell (z) = [\Wc_0] - [\Wc_{\ify}] =  
j_{0*} [\Xc] - 
j_{\ify *} [\Xc']\, . \leqno (14)
$$

Assume a class $\wt{\phi}_{\Xc_0} (E_{\Xc}^{\bu})$ 
satisfying properties i), ii), iii) has been defined. Let $p 
: \Wc \ra \Xc$ and $\pi : \Xc' \ra \Xc$
 be the projection maps. We get
$$
\leqalignno{
\ &p_* (i^* \, \phi ({E}_{\Wc}^{\bu}) \cdot \ell (z)) &(15) \cr
= \ &p_* (i^* \, i_* \, \phi_{\Wc_0} ({E}_{\Wc}^{\bu}) \cdot 
\ell (z)) \cr
= \ &p_* (\phi_{\Wc_0} ({E}_{\Wc}^{\bu}) \cdot i^* \, i_* \, 
\ell (z)) \cr
= \ &p_* (\phi_{\Wc_0} ({E}_{\Wc}^{\bu}) (j_{0*} [\Xc] - 
j_{\ify *} [\Xc'])) \cr
= \ & \phi_{\Xc_0} (j_0^* \, {E}_{\Wc}^{\bu}) - 
\pi_* (\phi_{\Xc'_0} (j_{\ify}^* \, {E}_{\Wc}^{\bu}) )\cr
= \ & \phi_{\Xc_0} ({E}_{\Xc}^{\bu}). \cr
}
$$
This proves the uniqueness of $\phi_{\Xc_0} (E_{\Xc}^{\bu})$. 
Note that this proof of uniqueness is the same as the one in 
[GS2] 1.3.2 for the archimedean analog. 

Conversely, if we define $\phi_{\Xc_0} (E_{\Xc}^{\bu})$
by formula (15) we can check properties i) to v) as in loc.cit..
Indeed, i) follows from the equalities

$$
\leqalignno{
\ &i_*p_* (i^* \, \phi ({E}_{\Wc}^{\bu}) \cdot \ell (z)) &(16) \cr
= \ &p_* ( \phi ({E}_{\Wc}^{\bu}) \cdot 
i_*\ell (z)) \cr
= \ &p_* ( \phi ({E}_{\Wc}^{\bu}) (j_{0*} [\Xc] - 
j_{\ify *} [\Xc'])) \cr
= \ & \phi (j_0^* \, {E}_{\Wc}^{\bu})\cr
= \ & \phi ({E}_{\Xc}^{\bu}). \cr
}
$$

Property ii) is clear. Under the assumption of iii),
there exists a split acyclic complex
 $({E}_{\Wc}^{\bu},d_{\Wc})$ extending
 $({E}_{\Xc}^{\bu},d_{\Xc})$
to $\Wc$, hence $\phi ({E}_{\Wc}^{\bu}) = 0$
and, by (15), $\phi_{\Xc_0} (E_{\Xc}^{\bu})$
vanishes. To prove (10), we note that
$$p_* (i^* \, ch ({E}_{\Wc}^{\bu}\ot p^*( F_{\Xc}) )
 \cdot \ell (z)) 
= {\rm ch}_{\Xc_0} (E_{\Xc}^{\bu} \ot F_{\Xc}).$$

It remains to prove the identities (8) and (9) for the behaviour of 
$\phi_{\Xc_0}$ in exact sequences. By deformation 
as in [Fu], proof of Proposition 18.1~b), 
or by iii) and the 
analog of [GS2] Prop.~1.3.4, we are reduced to the case where 
$(E_{\Xc}^{\bu} , d_{\Xc}) = (S_{\Xc}^{\bu} , d_{\Xc}) \op 
(Q_{\Xc}^{\bu} , d_{\Xc})$. But then we can copy 
 the proof of [GS2] Prop.~1.3.2. Namely let $\phi 
(S_{\Xc}^{\bu} , Q_{\Xc}^{\bu})$ be the class

$$
\phi (S_{\Xc}^{\bu} , Q_{\Xc}^{\bu})= {\rm ch}_{\Xc_0} 
(S_{\Xc}^{\bu}) + {\rm ch}_{\Xc_0} (Q_{\Xc}^{\bu}) $$
when $\phi = {\rm ch}$, and

$$
\phi (S_{\Xc}^{\bu} , Q_{\Xc}^{\bu}) = ({\rm Td} - 
1)_{\Xc_0} (S_{\Xc}^{\bu}) (i^* \, {\rm Td} (Q_{\Xc}^{\bu})) 
+ ({\rm Td} - 1)_{\Xc_0} (Q_{\Xc}^{\bu}) \, 
$$
when $\phi = {\rm Td} - 1 $.
This class is functorial and such that
$$
i_* \, \phi (S_{\Xc}^{\bu} , Q_{\Xc}^{\bu}) = \phi 
(S_{\Xc}^{\bu} \op Q_{\Xc}^{\bu}) \, .
$$
Furthermore it vanishes when there exist
acyclic complexes   $(S_{\Xc}^{\bu} , 
d_{\Xc})$ and $(Q_{\Xc}^{\bu} , d_{\Xc})$. Let 
$\Wc$ be a model of $X \build \ts_{K}^{} \Pb^1$
 together with maps $j_0 : 
\Xc \ra \Wc$, $j_{\ify} : \Xc' \ra \Wc$, $\pi : \Xc' \ra 
\Xc$ as above, where $\Wc$ maps to the closures of both
$
\vp_S (X \build \ts_{K}^{} \Pb^1)$ and 
$\vp_Q (X \build \ts_{K}^{} \Pb^1$, where
$$
\vp_S : X \build \ts_{K}^{} \Pb^1 \ra G \build \ts_{\Lb}^{} \Pb^1$$
and 
$$\vp_Q : X \build \ts_{K}^{} \Pb^1 \ra G' \build \ts_{\Lb}^{} \Pb^1
$$
are the maps defined as in [Fu] 18.1 from the complexes 
$(S^{\bu} , d)$ and $(Q^{\bu} , d)$ 
respectively. Let ${S}_{\Wc}^{\bu}$
and ${Q}_{\Wc}^{\bu}$ be the 
bundles on $\Wc$ defined as was ${E}_{\Wc}^{\bu}$
at the beginning of
this proof, and let 
$p : \Wc \ra \Xc$ be the projection map. We 
 get successively
$$
\eqalign{
\ &{\phi}_{\Xc_0} (S_{\Xc}^{\bu} \op Q_{\Xc}^{\bu}) \cr
= \ &p_* (i^* \, \phi ({S}_{\Wc}^{\bu} \op {Q}_{\Wc}^{\bu})
\cdot 
\ell (z)) \cr
= \ &p_* (i^* \, i_* \, \phi (S_{\Wc}^{\bu} ,Q_{\Wc}^{\bu})
\cdot \ell (z)) \cr
= \ &j_0^* \, \phi (S_{\Wc}^{\bu} , Q_{\Wc}^{\bu}) - 
\pi_* j_{\ify}^* \, \phi ({S}_{\Wc}^{\bu}, {Q}_{\Wc}^{\bu})\cr
= \ &\phi (S_{\Xc}^{\bu} , Q_{\Xc}^{\bu}) \, . \cr
}
$$
\hfill q.e.d.

\medskip

\noindent {\bf 2.4.} Let $\phi = {\rm ch}$ or ${\rm Td} - 
1$. It follows from Proposition~2 iv) that, if 
$(E_{\Xc}^{\bu} , d_{\Xc})$ and $(F_{\Xc}^{\bu} , d_{\Xc})$ 
are quasi-isomorphic complexes which are both acyclic on 
$X$,
$$
\phi_{\Xc_0} (E_{\Xc}^{\bu}) = \phi_{\Xc_0} (F_{\Xc}^{\bu}) 
\, .
$$
This implies that $\phi_{\Xc_0}$ defines a class 
$\phi_{\Xc_0} (K_{\Xc}) \in {\rm CH}_{\bu} (\Xc_0)_{\Qb}$ 
for any element $K_{\Xc}$ in $D_{\Xc_0}^b (\Xc)$ and, given 
a distinguished triangle
$$
K'_{\Xc} \ra K_{\Xc} \ra K''_{\Xc} \ra K'_{\Xc} \, [1]
$$
in $D_{\Xc_0}^b (\Xc)$, the formulae
$$
{\rm ch}_{\Xc_0} (K_{\Xc}) = {\rm ch}_{\Xc_0} (K'_{\Xc}) + 
{\rm ch}_{\Xc_0} (K''_{\Xc}) \leqno (17)
$$
and
$$
({\rm Td} - 1)_{\Xc_0} (K_{\Xc}) = ({\rm Td} - 1)_{\Xc_0} 
(K'_{\Xc}) \, {\rm Td} (K''_{\Xc}) + ({\rm Td} - 1)_{\Xc_0} 
(K''_{\Xc}) \leqno (18)
$$
hold, as well as
$$
i_* \, \phi_{\Xc_0} (K_{\Xc}^{\bu}) = \phi (K_{\Xc}^{\bu})\  .  
\leqno (19)
$$

Now let
$$
\Tc : K' \ra K \ra K'' \ra K' [1] \leqno (20)
$$
be a distinguished triangle in $D^b (X)$ and assume $K'$, 
$K$ and $K''$ are equipped with arbitrary metrics. We can 
attach to these data  classes $\wt{\rm ch} (\ov{\Tc})$, 
${\rm Td}_0 (\ov{\Tc})$ and $\wt{\rm Td} (\ov{\Tc})$ in $\wt 
A (X)_{\Qb}$ which are defined as follows. Choose an
exact sequence of complexes on $X$
$$
0 \ra E'^{\bu}\ra E^{\bu} \ra E''^{\bu} \ra 0 \leqno (21)
$$
representing $\Tc$. We view (21) as a double
complex where the second degree of $E'^{\bu}$ 
is  
zero, and we denote by $C^{\bu}$ the total 
complex of (21). The complex $C^{\bu}$
is acyclic and the metrics chosen on 
 $K'$, $K$ and $K''$ define a metric on
$C^{\bu}$, i.e. bundles $C^{n}_{\Xc} $
on some model $\Xc$  which restrict to 
$C^{n}$ on $X$ for all $n \geq 0$.
We then define
$$
\leqalignno{
\wt{\rm ch} (\ov{\Tc}) := & \ 
{\rm ch}_{\Xc_0} (C^{\bu}_{\Xc} ) &(22) \cr
{\rm Td}_0 (\ov{\Tc}) := & \ 
({\rm Td}-1)_{\Xc_0} (C^{\bu}_{\Xc} ) &(23) \cr
\wt{\rm Td} (\ov{\Tc}) := & \ {\rm Td}_0 (\ov{\Tc})  \, {\rm Td} 
(\ov K) \,  &(24) \cr
}
$$
in $\wt 
A (X)_{\Qb}$, where ${\rm Td} (\ov K) \in A_{\rm closed} (X)_{\Qb}$ 
is the Todd form of
$K$ with its chosen metric (cf.  \S 2.2).
These classes depend only on $\Tc$
and on the choice of metrics
on $K'$, $K$ and $K''$.
When $\phi = {\rm ch}$ or ${\rm Td}$, the following equation 
is satisfied:
$$
a (\wt{\phi} (\ov{\Tc})) = \wh{\phi} (\ov{K}' \op \ov{K}'') 
- \wh{\phi} (\ov K) \leqno (25)
$$
(this follows from (17)).

\smallskip

\smallskip

Let $h_0$ and $h_1$ be two metrics on a given $K \in D^b 
(X)$. Consider the triangle $\Tc$ as in (20) where $K'' = 0$, 
$K' = K$, and $K' \ra K$ is the identity. Let $K'$ be 
metrized by $h_0$ and $K$ by $h_1$. When $\phi = {\rm ch}$ 
or ${\rm Td}$ we define
$$
\wt{\phi} (h_0 , h_1) := \wt{\phi} (\ov{\Tc}) \, , \leqno (26)
$$
hence
$$
a (\wt{\phi} (h_0 , h_1)) = \wh{\phi} (K , h_0) - \wh{\phi} 
(K , h_1) \, . \leqno (27)
$$
It follows from (25) and (26) that, given three metrics $h_0$, 
$h_1$, $h_2$, we have
$$
\wt{\phi} (h_0 , h_1) = \wt{\phi} (h_0 , h_1) + \wt{\phi} 
(h_1 , h_2) \leqno (28)
$$
and that
$$
\wt{\phi} (h_0 , h_1) = - \wt{\phi} (h_1 , h_0) \, . \leqno 
(29)
$$

\smallskip

\noindent {\bf 2.5.} We also have the following

\medskip

\noindent {\bf Lemma 1.} {\it Let $\Tc$ be a triangle as in 
$(20)$. Let $h'_0$, $h_0$, $h''_0$ and $h'_1$, $h_1$, $h''_1$ be 
two sets of metrics on $K'$, $K$, $K''$ respectively. When 
$\phi = {\rm Td}$ or ${\rm ch}$, let $\wt{\phi} 
(\ov{\Tc}_0)$ and $\wt{\phi} (\ov{\Tc}_1)$ be the class defined 
above for each choice of metrics.}

\item{i)} {\it If $\phi = {\rm ch}$ we have}
$$
\wt{\rm ch} (\ov{\Tc}_0) - \wt{\rm ch} (\ov{\Tc}_1) = 
\wt{\rm ch} (h'_0 , h'_1) - \wt{\rm ch} (h_0 , h_1) + 
\wt{\rm ch} (h''_0 , h''_1) \, . \leqno (30)
$$

\item{ii)} {\it When $\phi = {\rm Td}$, $h_0 = h_1$ and 
$h''_0 = h''_1$, we have}
$$
\wt{\rm Td} (\ov{\Tc}_0) - \wt{\rm Td} (\ov{\Tc}_1) =  
\wt{\rm Td} (h'_0 , h'_1) \, {\rm Td} (\ov{K}'') \, . \leqno 
(31)
$$

\medskip

\noindent {\it Proof.} Let $\Xc$ be a model on which all the metrics
$h'_0$ ... $h''_1$ 
are defined, let 
$$
0 \ra E'^{\bu}\ra E^{\bu} \ra E''^{\bu} \ra 0 
$$
be an exact sequence representing
 $\Tc$ as in (21), and let
 $T_0$ (resp. 
$T_1$) be the associated total complex 
equipped with the metrics induced by
$h'_0$, $h_0$, $h''_0$  (resp. $h'_1$, $h_1$, $h''_1$). 
 Finally let $T$ be the total complex of the identity
map $T_0 \ra T_1$, equipped with the induced metric. 
 From the exact sequence
$$
0 \ra T_1 [1] \ra T \ra T_0 \ra 0 \leqno (32)
$$
 we get, using (8) and (11),
$$
{\rm ch}_{\Xc_0} (T) = {\rm ch}_{\Xc_0} (T_0) - {\rm 
ch}_{\Xc_0} (T_1) \, . \leqno (33)
$$
On the other hand, if $C'$ (resp. $C$, resp. $C''$)
is the cone of the identity map
$(E'^{\bu},h'_0) \ra (E'^{\bu},h'_1) $ 
(resp. $(E^{\bu},h_0) \ra (E^{\bu},h_1)$,
resp. $(E''^{\bu},h''_0) \ra (E''^{\bu},h''_1)$)
there are exact sequences
$$
0 \ra C'' [2] \ra T \ra T' \ra 0
$$
$$
0 \ra C[1] \ra T' \ra C' \ra 0
$$
from which it follows, by (8) and (11), that
$$
{\rm ch}_{\Xc_0} (T) = {\rm ch}_{\Xc_0} (C') - {\rm 
ch}_{\Xc_0} (C) + {\rm ch}_{\Xc_0} (C'') \, . \leqno (34)
$$
The equality (30) follows from (33) and (34).

\smallskip

On the other hand, we deduce from (32), (9) and (12) that
$$
({\rm Td} - 1)_{\Xc_0} (T) = ({\rm Td} - 1)_{\Xc_0} (T_0) \, 
{\rm Td} (T_1)^{-1} - ({\rm Td} - 1)_{\Xc_0} (T_1) \, {\rm 
Td} (T_1)^{-1} \, . \leqno (35)
$$
When $h_0 = h_1$ and $h''_0 = h''_1$, we can assume that $C$ 
and $C''$ are acyclic, therefore $T$ is quasi-isomorphic to 
$C'$ and
$$
({\rm Td} - 1)_{\Xc_0} (T) = ({\rm Td} - 1)_{\Xc_0} (C') \, 
. \leqno (36)
$$
If we multiply (35) and (36) by ${\rm Td} ( \ov{K}') \, {\rm 
Td} ( \ov{K}'')$ the identity (31) follows by (24). \hfill 
q.e.d.

\medskip

\noindent {\bf 2.6.} We  now define the (non-archimedean)
{\it arithmetic 
$K$-groups}. The group $\wh{K}_0 (X)$ is generated by 
triples $(E,h,\eta)$, where $E$ is a vector bundle on $X$, 
$h$ is a metric on $E$, and $\eta \in \wt{A} (X)_{\Qb}$. 
These generators are required to satisfy the following
relations. Let
$$
\Ec : 0 \ra S \ra E \ra Q \ra 0
$$
be any exact sequence of bundles on $X$ and suppose that 
$S$, $E$ and $Q$ are equipped with arbitrary metrics. Let 
$\wt{\rm ch} (\ov{\Ec}) \in \wt{A} (X)_{\Qb}$ be the 
corresponding secondary class, defined as in 2.4 (with $S$ 
in degree zero). Then, for any $\eta' \in  \wt{A} (X)_{\Qb}$ and 
$\eta'' \in  \wt{A} (X)_{\Qb}$, we have
$$
(\ov S , \eta') + (\ov Q , \eta'') = (\ov E , \eta' + \eta'' 
+ \wt{\rm ch} (\ov{\Ec})) \leqno (37)
$$
in $\wh{K}_0 (X)$.

\smallskip

If, in this definition, we allow $\eta$, $\eta'$, $\eta''$ 
to be any currents in $\wt{D} (X)_{\Qb}$, we get another 
group, denoted ${\build {K_0}_{}^{\vee}} (X)$, which clearly contains 
$\wh{K}_0 (X)$ as a subgroup.

\smallskip

There are maps
$$
\a : \wt{D} (X)_{\Qb} \ra {\build {K_0}_{}^{\vee}} (X)
$$
and
$$
{\rm ch} : {\build {K_0}_{}^{\vee}} (X) \ra D_{\rm closed} (X)_{\Qb}
$$
defined by $\a (\eta) = (0,\eta)$ and
$$
{\rm ch} (\ov E , \eta) = {\rm ch} (\ov E) + dd^c \, \eta \, 
. \leqno (38)
$$
Note the equalities
$$
{\rm ch} \circ \a = \om \circ a = dd^c \, .
$$
Clearly, $x \in {\build {K_0}_{}^{\vee}} (X)$ lies in $\wh{K}_0 (X)$ if 
and only if ${\rm ch} (x)$ lies in $A_{\rm closed } 
(X)_{\Qb}$. From Proposition~1 i), we get that $\a (\eta)$ 
lies in $\wh{K}_0 (X)$ if and only if $\eta$ lies in $\wt{A} 
(X)_{\Qb}$.

\smallskip

From Proposition~2 it follows that any element $K$ in $D^b 
(X)$, when equipped with a metric as in 2.1, has a class 
$[\ov K]$ in $\wh{K}_0 (X)$ and that, if $K$ is acyclic,
$$
[\ov K] = \a \, (\wt{\rm ch} \, (\ov K)) \, .
$$

Since forms and metrized vector bundles are contravariant, 
any map of varieties $f : Y \ra X$ induces a pull-back 
morphism
$$
f^* : \wh{K}_0 (X) \ra \wh{K}_0 (Y) \, .
$$

The formula
$$
(\ov E , \eta) \cdot (\ov F , \xi) = (\ov E \ot \ov F , {\rm 
ch} (\ov E) \, \xi + \eta \, {\rm ch} (\ov F) + \eta \, dd^c 
\, \xi)
$$
defines both a ring structure on $\wh{K}_0 (X)$ and a module 
structure of ${\build {K_0}_{}^{\vee}} (X)$ over $\wh{K}_0 (X)$ (that 
this map is compatible with (37) is a consequence of 
Proposition~2 v)). Note that $f^* (xy) = f^* (x) \, f^* 
(y)$ when $x$ and $y$ lie in $\wh{K}_0 (X)$, and that
$$
x \, \a (\eta) = \a ({\rm ch} (x) \, \eta)\,  \leqno (39)$$
when $\eta \in \wt D (X)_{\Qb}$ and  
$ x \in \wh{K}_0 (X)$.

There is a Chern character map
$$
{\build {\rm ch}_{}^{\vee}} : {\build {K_0}_{}^{\vee}} (X) 
\ra {\build {{\rm CH}^{\bu}}_{}^{\vee}} 
(X)_{\Qb} \, ,
$$
defined by mapping $(\ov E , \eta)$ to the class of $\wh{\rm 
ch} (\ov E) + a (\eta)$. It induces a ring homomorphism
$$
\wh{\rm ch} : \wh{K}_0 (X) \ra \wh{\rm CH}^{\bu} (X)_{\Qb}
$$
which commutes with pull-backs. Note that
$$
\om \circ {\build {\rm ch}_{}^{\vee}} = {\rm ch} \qquad \hbox{and} \qquad 
{\build {\rm ch}_{}^{\vee}} \circ \a = a \, . \leqno (40)
$$

Our main goal will be to define push-forward morphisms 
$f_*$ on arithmetic $K$-groups, satisfying a Riemann-Roch 
formula. For that purpose we need preliminaries on tangent 
complexes.

\bigskip

\noindent {\bf 3. Secondary Todd classes of tangent 
complexes}

\medskip

\noindent {\bf 3.1.} Given any map of varieties $\vp : X \ra 
Y$ (resp. any map of models $f : \Xc \ra \Yc$) we can attach 
to it a tangent complex $T\vp \in D^b (X)$ (resp. $Tf \in D^b 
(\Xc)$). A representative of $T\vp$ is
 the complex of vector bundles $TX \ra \vp^*(TY)$
with $TX$ in degree zero.
Similarly, when both models are smooth over some
base, $Tf$ is represented by $T\Xc \ra f^*(T\Yc)$;
in general it is defined by the  construction
 dual to the one of the cotangent complex
in [SGA6] Exp. VIII, \S 2.

\smallskip

Given any map of models $f : \Xc \ra \Yc$, we shall denote 
by
$$
{\rm Td} (f) \in {\rm CH}^{\bu} (\Xc_0)_{\Qb} \sbs A_{\rm 
closed}^{\bu} (X)_{\Qb}
$$
the Todd class $i^* \, {\rm Td} (Tf)$. When $\pi : \Xc' \ra 
\Xc$ is a morphism between models of a given variety
$X$, the tangent complex 
$T\pi$ is acyclic on $X$. Therefore it defines a class
$$
{\rm Td}_0 (\pi) := ({\rm Td}-1)_{\Xc'_0} (T\pi) \in {\rm 
CH}_{\bu} (\Xc'_0)_{\Qb} \sbs \wt{A} (X)_{\Qb} \, .
$$
Finally, when $f : \Xc \ra \Yc$ and $f' : \Xc' \ra \Yc'$ are 
two maps of models which induce the same map $\vp : X \ra Y$
of varieties, 
we let
$$
\wt{\rm Td} (f,f') := \wt{\rm Td} (h_0 , h_1) \in \wt A 
(X)_{\Qb} \, , \leqno (41)
$$
where $h_0$ (resp. $h_1$) is the metric induced by $Tf$ 
(resp. $Tf'$) on $T{\vp}$.

 When $f : \Xc \ra \Yc$ and $g : 
\Yc \ra \Zc$ are two map of models, there is a distinguished 
triangle in $D^b (\Xc)$:
$$
Tf \ra T(gf) \ra f^* \, Tg \ra Tf [1] \, . \leqno (42)
$$
It follows that
$$
{\rm Td} (gf) = {\rm Td} (f) \, f^* \, {\rm Td} (g) \, . \leqno 
(43)
$$
This also inplies that, if $\pi : \Xc' \ra \Xc$ is a 
morphism between two models of the variety
$X$ and $f : \Xc \ra \Yc$ is any 
map of models, the cone
$$
T (f \pi) \ra \pi^* \, Tf
$$
is isomorphic to $T \, \pi \, [1]$ and therefore, 
from (24) and from (41), the following identity holds
true in $\wt A 
(X)_{\Qb}$:
$$
{\rm Td}_0 (\pi) \, \pi^* ({\rm Td} (f)) = \wt{\rm Td} (f 
\pi , f) \, . \leqno (44)
$$
Similarly, if $g : \Yc \ra \Xc'$ is any map of models, the 
identity
$$
g^* ({\rm Td}_0 (\pi)) \cdot {\rm Td} (g) = \wt{Td} (\pi g , 
g) \leqno (45)
$$
holds in $\wt A (Y)_{\Qb}$.

\smallskip

Finally, when $\pi : \Xc' \ra \Xc$ and $\rho : \Xc'' \ra 
\Xc'$ are two maps between models of a given variety $X$, we 
deduce from (42) and (18) that
$$
{\rm Td}_0 (\pi \rho) = {\rm Td}_0 (\rho) \, \rho^* \, {\rm 
Td} (\pi) + \rho^* \, {\rm Td}_0 (\pi) \, . \leqno (46)
$$

\smallskip

\noindent {\bf 3.2. Proposition 3.} {\it Let $\pi : \Xc' \ra 
\Xc$ be a  map between two models of $X$.}

\item{i)} {\it For any
$K \in D^b (\Xc)$ the adjunction map}
$
K \ra R\pi_* \, L\pi^* \, K \,  
$
{\it is an isomorphism.}

\item{ii)} {\it
Let ${\rm 
Td}_0 (\pi) = ({\rm Td}-1)_{\Xc'_0} (T\pi) \in {\rm CH}_{\bu} 
(\Xc'_0)_{\Qb}$. Then, in ${\rm CH}_{\bu} (\Xc_0)_{\Qb}$, we 
have}
$$
\pi_* \, {\rm Td}_0 (\pi) = 0 \, . \leqno (47)
$$

\smallskip

\noindent {\it Proof.} To prove i), note that $R\pi_* \, L\pi^* \, K$
is the derived tensor product of $K$ with $ R \pi_* \,
\Oc_{\Xc'} $, hence it is enough to consider the case where
$K= \Oc_{\Xc} $. If we assume furthermore that
$\pi$ is the blow up of a closed regular subscheme in
${\Xc} $, the assertion follows from [SGA6] Exp.~VII, Lemme~3.5, 
p.~441. Therefore i) is true when $\pi$ is a good map of models in
the sense of \S 1.1. Using the axiom (M2) of loc. cit., 
we may find two models ${\Xc}_1, {\Xc}_2 $, and maps of models
${\Xc}_2 \ra {\Xc}_1 $ and ${\Xc}_1 \ra {\Xc'} $
so that the composite maps ${\Xc}_2 \ra {\Xc'} $
and ${\Xc}_1 \ra {\Xc} $ are good. Let $\pi_1 :{\Xc}_1 \ra {\Xc} $
and $\pi_2 :{\Xc}_2 \ra {\Xc} $ be the two obvious maps of models.
We get  a sequence of morphisms in $D^b (\Xc)$ :
$$\Oc_{\Xc} \ra  R \pi_* \, \Oc_{\Xc'} \ra  R \pi_{1*} \, \Oc_{{\Xc}_1}
\ra  R \pi_{2*} \, \Oc_{{\Xc}_2 } \, .$$
The composite of the first two maps is the isomorphism
$\Oc_{\Xc} \ra  R \pi_{1*} \, \Oc_{{\Xc}_1}$ (since $\pi_1$ is good)
and the composite of the last two maps
is also an isomorphism $R \pi_* \, \Oc_{\Xc'}
\ra  R \pi_{2*} \, \Oc_{{\Xc}_2}$ since ${\Xc}_2 \ra {\Xc'} $
is good. It follows that all morphisms in the sequence above
are isomorphisms. This proves i).

To prove ii), we apply the refined Riemann-Roch formula conjectured by
T.Saito [S] p.163, and proved by J.Franke [Fr] \S3.3. 
Consider the statement in [S] loc.cit. when 
$Y = T = {\Xc}$, $R= {\Xc_0}$,
 $Z = {\Xc'}$, $ h = id_{\Xc}$, $ \pi = g$, and 
$F = \Oc_{\Xc}$. By i) above, the canonical map 
$$
Rh_* \, F \ra Rg_* \, L\pi^* \, F \,  
$$
 is an isomorphism in that case, so the left hand side 
 of [S], loc. cit., vanishes. On the other hand, 
 the right hand side is precisely
 $\pi_* \, {\rm Td}_0 (\pi)$. 

 \hfill q.e.d.

\smallskip

\noindent {\bf 3.3.}
When $\pi$ is  a good map of models, one can also prove
Proposition 3 ii) directly. Indeed, by definition,
 $\pi$ is then the composite of 
a sequence of blow ups with 
integral regular centers meeting the reduced 
special fiber normally. By (46) and the projection formula, it 
is enough to check (47) when $\pi$ is one of these good blow 
ups. Let $\Yc$ be the center of this blow up. We let $j_0 : 
\Yc \ra \Xc_0$ and $j : \Yc \ra \Xc$ be the obvious 
inclusions and denote by $\Yc' = \pi^{-1} (\Yc)$ the 
exceptional divisor of $\pi$. If $N$ is the normal bundle of 
$\Yc$ in $\Xc$, we know that $\Yc' = \Pb (N)$. Let 
 $j' : \Yc' \ra \Xc'$ be the 
inclusion. According to [F], Lemma 15.4 (iv),
 the tangent complex $T\pi$, when shifted by one,
is canonically isomorphic the direct image
 $j'_* (F)$ of the universal quotient bundle $F$
 on $\Pb (N)$.
Therefore, if we apply the Grothendieck-Riemann-Roch theorem 
with supports
to $j'_*$, we see that ${\rm Td}_0 (\pi)$ is $j'_* (\tau)$, 
where $\tau$ is a universal polynomial in the Chern classes 
of $p^* N$ and the Chern class of the canonical line bundle 
$\Oc (1)$ on $\Pb (N)$. It follows that $\pi_* \, {\rm Td}_0 
(\pi) = j_{0*} (\s)$, where $\s$ is a universal polynomial 
$R(c_1 (N), \ldots , c_r (N))$ in the Chern classes of $N$, 
where $R$ depends only on the rank $r$ of $N$.

\smallskip

On the other hand, if we apply the Grothendieck-Riemann-Roch 
theorem to $\pi$ and $\Oc_{\Xc'}$, since 
$R \pi_* \, 
\Oc_{\Xc'} = \Oc_{\Xc}$, we 
obtain 
$$
1 = {\rm ch} (\Oc_{\Xc}) = \pi_* \, {\rm Td} (\pi) \, ,
$$
therefore, by (19),
$$
i_* \, \pi_* \, {\rm Td}_0 (\pi) = \pi_* ({\rm Td} (\pi) - 
1) = 0 \, .
$$
In particular Proposition 3 ii) holds as soon as $i_*$ is injective, and 
$\s = 0$ when $j_*$ is injective. The polynomial $R$ is the 
same for any regular closed immersion $j$ ($\Xc$ need not 
defined over $\Lb$), for instance the standard section 
$j : 
\Yc \ra \Pb (N \op 1)$ of the completed projective bundle of 
$N$, for which $j_*$ is injective. Therefore this universal 
polynomial $R$ must vanish. Hence we always have $\s = 0$, 
and Proposition~3 holds.

\smallskip

\noindent {\bf 3.4.} When $\Lb$ is a 
localization of an algebra of finite type over a field $k$
of 
characteristic zero, the general case of Proposition 3 ii) would 
follow from the case of good maps of models and the strong
factorization conjecture for birational maps made in
[AKMW] (0.2.1) and [W].
Indeed, by a standard inductive limit argument,
one is reduced to proving 
(47) when $\pi : \Xc' \ra 
\Xc$ 
 is a birational map between 
two smooth projective varieties over $k$
which is the identity outside the closed
subset $\Xc_0 \subset \Xc$.
But, according to the strong factorization
conjecture, there exist a 
smooth projective variety $\Xc''$
and two maps $f: \Xc'' \ra \Xc$ and
$g: \Xc'' \ra \Xc'$ which are  compositions of a sequence
of blow ups with smooth centers contained 
in the inverse image of $\Xc_0$, so that $f= \pi \circ g$.
The identity (47) is true for $f$ and $g$ by \S 3.3, and 
(46) implies it for $\pi$.

 We have not been able to prove (47)
by using only \S 3.3 and the weak factorization 
theorem of [W] and [AKMW].

\bigskip

\noindent {\bf 4. Direct images}

\smallskip

\noindent {\bf 4.1.} Let $\vp : X \ra Y$ be a map of 
varieties and consider a commutative diagram of maps of 
models
$$
\diagram{
\Xc' &\hfl{\pi}{} &\Xc \cr
\vfl{f'}{} &&\vfl{}{f} \cr
\Yc' &\hfl{p}{} &\Yc \cr
} \leqno (48)
$$
where $\Xc$ and $\Xc'$ (resp. $\Yc$ and $\Yc'$) are models 
of $X$ (resp. $Y$) and both $f$ and $f'$ induce the same map 
$\vp$ from $X$ to $Y$.  
Let $E_{\Xc}$ be a vector bundle on $\Xc$. The elements 
$Lp^* \, Rf_* \, E_{\Xc}$ and $Rf'_* \, L\pi^* \, E_{\Xc}$ 
in $D^b (\Yc')$ both restrict to $R \vp_* \, E$ on $Y$, when 
$E$ is the restriction of $E_{\Xc}$ to $X$. Furthermore 
there is a canonical morphism
$$
Lp^* \, Rf_* \, E_{\Xc} \ra Rf'_* \, L\pi^* \, E_{\Xc} \, . 
\leqno (49)
$$
Indeed, let $\Xc_1$ be the fiber product of $\Xc$ and $\Yc'$ 
over $\Yc$, and $\ve : \Xc' \ra \Xc_1$, $\pi_1 : \Xc_1 \ra 
\Xc$ and $f_1 : \Xc_1 \ra \Yc'$ the obvious morphisms. By 
adjunction, as in \S 3.3.,
 there is a morphism of functors ${\rm id} \ra R 
\ve_* \, L \ve^*$ in the derived category of perfect 
complexes on $\Xc_1$. This gives a map
$$
Rf_{1*} \, L\pi_1^* \, E_{\Xc} \ra Rf'_* \, L\pi^* \, 
E_{\Xc} = Rf_{1*} \, R\ve_* \, L\ve^* \, L\pi_1^* \, E_{\Xc} 
\, .
$$
The map (49) is the composite of this map with the base change
morphism
$$
Lp^* \, Rf_* \, E_{\Xc} \ra Rf_{1*} \, L\pi_1^* \, E_{\Xc}
$$
([SGA4] XVII 4.1.4). Let $K_{\Yc'} \in D_{\Yc'_0}^b (\Yc')$ 
be the cone of the map (49). We define
$$
\leqalignno{
\t_{\Yc'} := \ & {\rm ch}_{\Yc'_0} (K_{\Yc'}) \cdot {\rm Td} 
(p) &(50) \cr
+ \ &p^* ({\rm ch} (Rf_* \, E_{\Xc})) \, {\rm Td}_0 (p) \cr
- \ &f'_* [\pi^* ({\rm ch} (E_{\Xc}) \, {\rm Td} (f)) \cdot 
{\rm Td}_0 (\pi)] \cr
}
$$
in ${\rm CH}_{\bu} (\Yc'_0)_{\Qb}$.

\smallskip

Note that, when $\Yc' = \Yc$, we have
$$
\t_{\Yc} = 0 \, . \leqno (51)
$$

\smallskip

\bigskip

\noindent {\bf 4. Direct images}

\smallskip

\noindent {\bf 4.1.} Let $\vp : X \ra Y$ be a map of 
varieties and consider a commutative diagram of maps of 
models
$$
\diagram{
\Xc' &\hfl{\pi}{} &\Xc \cr
\vfl{f'}{} &&\vfl{}{f} \cr
\Yc' &\hfl{p}{} &\Yc \cr
} \leqno (48)
$$
where $\Xc$ and $\Xc'$ (resp. $\Yc$ and $\Yc'$) are models 
of $X$ (resp. $Y$) and both $f$ and $f'$ induce the same map 
$\vp$ from $X$ to $Y$.  
Let $E_{\Xc}$ be a vector bundle on $\Xc$. The elements 
$Lp^* \, Rf_* \, E_{\Xc}$ and $Rf'_* \, L\pi^* \, E_{\Xc}$ 
in $D^b (\Yc')$ both restrict to $R \vp_* \, E$ on $Y$, when 
$E$ is the restriction of $E_{\Xc}$ to $X$. Furthermore 
there is a canonical morphism
$$
Lp^* \, Rf_* \, E_{\Xc} \ra Rf'_* \, L\pi^* \, E_{\Xc} \, . 
\leqno (49)
$$
Indeed, let $\Xc_1$ be the fiber product of $\Xc$ and $\Yc'$ 
over $\Yc$, and $\ve : \Xc' \ra \Xc_1$, $\pi_1 : \Xc_1 \ra 
\Xc$ and $f_1 : \Xc_1 \ra \Yc'$ the obvious morphisms. By 
adjunction, as in \S 3.3.,
 there is a morphism of functors ${\rm id} \ra R 
\ve_* \, L \ve^*$ in the derived category of perfect 
complexes on $\Xc_1$. This gives a map
$$
Rf_{1*} \, L\pi_1^* \, E_{\Xc} \ra Rf'_* \, L\pi^* \, 
E_{\Xc} = Rf_{1*} \, R\ve_* \, L\ve^* \, L\pi_1^* \, E_{\Xc} 
\, .
$$
The map (49) is the composite of this map with the base change
morphism
$$
Lp^* \, Rf_* \, E_{\Xc} \ra Rf_{1*} \, L\pi_1^* \, E_{\Xc}
$$
([SGA4] XVII 4.1.4). Let $K_{\Yc'} \in D_{\Yc'_0}^b (\Yc')$ 
be the cone of the map (49). We define
$$
\leqalignno{
\t_{\Yc'} := \ & {\rm ch}_{\Yc'_0} (K_{\Yc'}) \cdot {\rm Td} 
(p) &(50) \cr
+ \ &p^* ({\rm ch} (Rf_* \, E_{\Xc})) \, {\rm Td}_0 (p) \cr
- \ &f'_* [\pi^* ({\rm ch} (E_{\Xc}) \, {\rm Td} (f)) \cdot 
{\rm Td}_0 (\pi)] \cr
}
$$
in ${\rm CH}_{\bu} (\Yc'_0)_{\Qb}$.

\smallskip

Note that, when $\Yc' = \Yc$, we have
$$
\t_{\Yc} = 0 \, . \leqno (51)
$$

\smallskip

\noindent {\bf Proposition 4.} 

\item{i)} {\it When $f : \Xc \ra \Yc$ and $E_{\Xc}$ are 
fixed, and when $f' : \Xc' \ra \Yc'$ varies, the classes 
$\t_{\Yc'}$ are the components of a unique current $\t \in 
\wt D (Y)_{\Qb}$.}

\item{ii)} {\it The following identity holds in ${\rm CH} 
(\Yc')_{\Qb}$:}
$$
i_* \, \t_{\Yc'} = f'_* \, \pi^* ({\rm ch} (E_{\Xc}) \, {\rm 
Td} (f)) - p^* \, {\rm ch} (Rf_* \, E_{\Xc}) \, . \leqno (52)
$$

\smallskip

\noindent {\it Proof.} By a remark made in \S~1.2, it 
follows from (M2) that, to check i), it is enough to prove 
the following. Consider a commutative diagram
$$
\diagram{
\Xc'' &\hfl{\rho}{} &\Xc' &\hfl{\pi}{} &\Xc \cr
\vfl{f''}{} &&\vfl{f'}{} &&\vfl{f}{} \cr
\Yc'' &\hfl{}{r} &\Yc' &\hfl{}{p} &\Yc \, , \cr
} \leqno (53)
$$
where $\rho$ and $\pi$ (resp. $r$ and $p$) are  maps of 
models of $X$ (resp. $Y$) and let $E_{\Xc}$ be a vector 
bundle on $\Xc$. Then
$$
r_* \, \t_{\Yc''} = \t_{\Yc'} \, , \leqno (54)
$$
where $\t_{\Yc'}$ is defined using $\pi$ and $p$, while 
$\t_{\Yc''}$ is defined with $\s = \pi \rho$ and $s = pr$. 

Let $K_{\Yc''}$ be the cone of the map
 $$Ls^* \, Rf_* \, E_{\Xc} 
\ra Rf''_* \, L\s^* \, E_{\Xc}.$$
 Note that
$$
Rr_* \, Ls^* = Rr_* \, Lr^* \, Lp^* = Lp^* \, . 
$$
Indeed, we know from Proposition 3 i) that $Rr^* \, Lr^* = {\rm id}$.
 Therefore
$$
Rr_* \, Ls^* \, Rf_* \, E_{\Xc} = Lp^* \, Rf_* \, E_{\Xc} \, 
. \leqno (55)
$$

Similarly $R\rho_* \, L\rho^* = {\rm id}$, hence
$$
Rr_* \, Rf''_* \, L\s^* \, E_{\Xc} = Rf'_* \, R\rho_* \, 
L\s^* \, E_{\Xc} = Rf'_* \, L\pi^* \, E_{\Xc} \, . \leqno (56)
$$
From (55) and (56) it follows that
$$
Rr_* \, K_{\Yc''} = K_{\Yc'}
$$
and the Riemann-Roch-Grothendieck theorem with supports 
([Fu] Th. 18.2)
implies that, since ${\rm Td} (s) = {\rm Td} (r) \, r^* \, {\rm Td} (p)$,
$$
r_* ({\rm ch}_{\Yc''_0} (K_{\Yc''}) \, {\rm Td} (s)) = {\rm 
ch}_{\Yc'_0} (K_{\Yc'}) \, {\rm Td} (p) \, . \leqno (57)
$$

Let us  look at the other summands of $\t_{\Yc''}$
(see (50)). If we apply 
$r_*$ to the identity
$$
{\rm Td}_0 (s) = {\rm Td}_0 (r) \, r^* \, {\rm Td} (p) + r^* 
\, {\rm Td}_0 (p)
$$
(see (46)) we obtain, by Proposition~3 ii) for $r$,
$$
r_* \, {\rm Td}_0 (s) = {\rm Td}_0 (p) \, .
$$
It follows that
$$
\leqalignno{
\ &r_* (s^* ({\rm ch} \, Rf_* \, E_{\Xc}) \, {\rm Td}_0 (s)) 
&(58) \cr
= \ &p^* \, {\rm ch} (Rf_* \, E_{\Xc}) \, r_* \, {\rm Td}_0 
(s) \cr
= \ &p^* \, {\rm ch} (Rf_* \, E_{\Xc}) \, {\rm Td}_0 (p) \, 
. \cr
}
$$
Similarly $r_* \, f''_* = f'_* \, \rho_*$ and $\rho_* \, 
{\rm Td}_0 (\s) = {\rm Td}_0 (\pi)$, therefore
$$
\leqalignno{
\ &r_* \, f''_* (\s^* ({\rm ch} (E_{\Xc}) \, {\rm Td} (f)) 
\, {\rm Td}_0 (\s)) &(59) \cr
= \ &f'_* \, \rho_* (\s^* ({\rm ch} (E_{\Xc}) \, {\rm Td} 
(f)) \, {\rm Td}_0 (\s)) \cr
= \ &f'_* (\pi^* ({\rm ch} (E_{\Xc}) \, {\rm Td} (f)) \, 
\pi_* \, {\rm Td}_0 (\s)) \cr
= \ &f'_* (\pi^* ({\rm ch} (E_{\Xc}) \, {\rm Td} (f)) \, 
{\rm Td}_0 (\pi)) \, . \cr
}
$$
From (57), (58) and (59) we conclude that $r_* \, \t_{\Yc''} = 
\t_{\Yc'}$, as was to be shown.

\smallskip

To prove ii) we use (17) and (50) to get
$$
\leqalignno{
i_* \, \t_{\Yc'} = \ &({\rm ch} (f'_* \, \pi^* \, E_{\Xc}) - 
{\rm ch} (p^* \, f_* \, E_{\Xc})) \, {\rm Td} (p) &(60) \cr
\ &+ p^* ({\rm ch} (f_* \, E_{\Xc})) \, ({\rm Td} (p) - 1) 
\cr
\ &- f'_* [\pi^* ({\rm ch} (E_{\Xc}) \, {\rm Td} (f)) \, 
({\rm Td} (\pi) - 1)] \cr
= \ &{\rm ch} (f'_* \, \pi^* \, E_{\Xc}) \, {\rm Td} (p) \cr
\ & - p^* \, {\rm ch} (f_* \, E_{\Xc}) \cr
\ &+ f'_* \, \pi^* ({\rm ch} (E_{\Xc}) \, {\rm Td} (f)) \cr
\ &- f'_* [\pi^* ({\rm ch} (E_{\Xc}) \, {\rm Td} (f)) \, 
{\rm Td} (\pi)] \, , \cr
}
$$
where we wrote $f_*$ instead of $Rf_*$, $p^*$ instead of 
$Lp^*$ etc. The relative Todd class of $f \circ \pi = p 
\circ f'$ is
$$
\pi^* ({\rm Td} (f)) \, {\rm Td} (\pi) = {\rm Td} (f') \, 
f'^* \, ({\rm Td} (p)) \, ,
$$
therefore the Riemann-Roch-Grothendieck theorem, when 
applied to $f'$, gives
$$
\leqalignno{
\ &f'_* [\pi^* ({\rm ch} (E_{\Xc}) \, {\rm Td} (f)) \, {\rm 
Td} (\pi)] &(61) \cr
= \ &f'_* [\pi^* ({\rm ch} (E_{\Xc})) \, {\rm Td} (f') \, 
f'^* \, ({\rm Td} (p))] \cr
= \ &{\rm ch} (f'_* \, \pi^* \, E_{\Xc}) \, {\rm Td} (p) \, 
. \cr
}
$$

Combining (60) and (61) we get (52). \hfill q.e.d.

\medskip

Let $\ov E$ (resp. $\ov{R\vp_* \, E}$) be the bundle $E$ 
(resp. $R\vp_* \, E$) equipped with the metric defined by 
$E_{\Xc}$ (resp. $Rf_* \, E_{\Xc}$).

\medskip

\noindent {\bf Corollary 2.} {\it The following identity of 
currents hold in $D_{\rm closed} (Y)_{\Qb}$:}
$$
dd^c \, \t = \vp_* ({\rm ch} (\ov E) \, {\rm Td} (f)) - {\rm 
ch} (\ov{R\vp_* \, E}) \, . \leqno (62)
$$

\smallskip

This corollary indicates that $\t$ plays the role of the
{\it 
higher analytic torsion} in Arakelov geometry ([GS3] or 
[BK]). Note that, when $f$ is flat, it follows from (), 
Prop.~1 i) and [BGS] Th.~4.1.1, that $\t$ lies in $\wt A 
(X)_{\Qb}$.

\medskip

\noindent {\bf 4.2. Theorem 1.} {\it Let $f : X \ra Y$ be a 
map of varieties. Choose a metric $h_f$ on the tangent 
complex $Tf$.
There exists a unique direct image morphism
$$
f_* : {\build {K_0}_{}^{\vee}} (X) \ra {\build {K_0}_{}^{\vee}} (Y)
$$
such that}

\item{i)} When $x = \a 
(\eta)$ with $\eta \in \wt D (X)_{\Qb}$, the following formula holds:
$$
f_* (\a (\eta)) = \a (f_* (\eta \, {\rm Td} (\ov{Tf}))) \, . 
\leqno (63)
$$

\item{ii)} {\it Assume there is a map of models
$$
\ud{f} : \Xc \ra \Yc
$$
such that $h_f$ is defined by $T\ud f$, and that $x \in 
{\build {K_0}_{}^{\vee}} (X)$ is the class of $(E,E_{\Xc},0)$, where 
$E_{\Xc}$ is a bundle on $\Xc$ with restriction $E$ to $X$. 
Let $\t$ be the current defined in Proposition~4, i). Then 
$f_* (x)$ is the class of $(Rf_* \, E, R\ud{f}_* \, E_{\Xc} 
, \t)$ in ${\build {K_0}_{}^{\vee}} (Y)$.}

\item{iii)} {\it Suppose we choose two different metrics 
$h_f$ and $h'_f$ on $Tf$ and let $f_*$, $f'_*$ be the 
corresponding direct image morphisms. Then, for any $x \in 
{\build {K_0}_{}^{\vee}} (X)$, the following identity holds in 
${\build {K_0}_{}^{\vee}} (Y)$:}
$$
f_* (x) - f'_* (x) = \a (f_* ({\rm ch} (x) \, \wt{\rm Td} 
(h_{f} , h'_f)) \, . \leqno (64)
$$

{\it Furthermore, for any $x$ in ${\build {K_0}_{}^{\vee}} (X)$ the 
following Riemann-Roch identity holds:}
$$
{\build {\rm ch}_{}^{\vee}} (f_* (x)) = f_* ({\build {\rm ch}_{}^{\vee}} (x) \, 
\wh{\rm Td} (\ov{Tf})) \, . \leqno (65)
$$
\smallskip

\noindent {\bf 4.3.} To prove uniqueness in Theorem~1, first 
notice that the identity (63) fixes $f_*$ on the image of $\a$.

\smallskip

Next, if $h$ and $h'$ are two metrics on a vector bundle $E$ 
over $X$, the relation (37) in ${\build {K_0}_{}^{\vee}} (X)$ together with 
(63) imply that if $x$ is the class of $(E,h,0)$ in 
${\build {K_0}_{}^{\vee}} (X)$ and $x'$ is the class of $(E,h',0)$, we 
must have
$$
f_* (x) - f_* (x') = \a (f_* (\wt{\rm ch} (h,h') \, {\rm Td} 
(\ov{Tf}))) \, . \leqno (66)
$$

On the other hand, given any map $f : X \ra Y$ of varieties 
and any vector bundle $E$ on $X$, we may find a map of 
models $\ud f : \Xc \ra \Yc$ inducing $f$ on $X$,
and a bundle 
$E_{\Xc}$ on $\Xc$ inducing $E$ on $X$ [RG]. If $Tf$ is 
equipped with the metric defined by $T \ud f$, Theorem ii)
will then specify the value of $f_* (x)$, where $x$ is the 
class of $(E,E_{\Xc},0)$. This, together with (63) and the 
anomaly formulae (64) and (66), proves the uniqueness of $f_*$.

\smallskip

To prove the existence of $f_*$ we have to show that the 
formula (63) and ii) are consistent with the anomaly formulae 
(64) and (66). This boils down to the following two facts. 
First, let (53) be the diagram considered in the proof of 
Proposition~4, and let $E_{\Xc}$ be a vector bundle on $\Xc$. Let 
$x \in {\build {K_0}_{}^{\vee}} (X)$ be the class of 
$(E,E_{\Xc},0)$ and let 
$\vp : X \ra Y$ be the map of varieties induced by $f$, $f'$ 
and $f''$. We want to compare the direct images of $x$ when 
the metric on $T\vp$ is defined by $Tf$ or $Tf'$. Let $\t$ 
be the current defined by $E_{\Xc}$ and the metric $Tf$, and 
$\t'$ the current defined by $\pi^* E_{\Xc}$ and $Tf'$ 
(Proposition~4). Let $h_0$ be the metric $Rf_* \, E_{\Xc}$ 
on $R\vp_* \, E$ and $h_1$ be the metric $Rf'_* \, L\pi^* \, 
E_{\Xc}$ on the same complex. Combining (37) and (64), the 
identity
$$
\t - \t' + \wt{\rm ch} (h_0 , h_1) = \vp_* ({\rm ch} (x) \, 
\wt{\rm Td} (Tf , Tf')) \, 
$$
must be true, at least after applying $\a$ to it.
To prove this equality, writing $K$ for the cone of $Lp^* \, Rf_* 
\, E_{\Xc} \ra Rf'_* \, L\pi^* \, E_{\Xc}$, all we need to 
check is the following identity in ${\rm CH}_{\bu} 
(\Yc''_0)_{\Qb}$:
$$
\t_{\Yc''} - \t'_{\Yc''} - r^* \, {\rm ch}_{\Yc''_0} (K) = 
f''_* ({\rm ch} (\s^* \, E_{\Xc}) \, \rho^* \, \wt{\rm Td} 
(f',f)_{\Yc'_0}) \, . \leqno (67)
$$
Conversely, by the cofinality axiom (M2),
 this identity will allow us to 
define $f_* (x)$ for any choice of metrics on $Tf$, and (64) 
will hold always.

\smallskip

Second, to check that $f_*$ defines a map on ${\build {K_0}_{}^{\vee}} (X)$, 
consider the diagram (48) in 4.1, let $S_{\Xc}$, 
$E_{\Xc}$, $Q_{\Xc}$ be bundles on $\Xc$, which restrict to 
$S$, $E$, $Q$ respectively on $X$. Assume there is a complex
$$
\Ec_{\Xc} : 0 \ra S_{\Xc} \ra E_{\Xc} \ra Q_{\Xc} \ra 0
$$
on $\Xc$ which restricts to an exact complex on $X$. 
Consider the associated Chern character with supports
$$
{\rm ch}_{\Xc_0} (\Ec_{\Xc}) \in {\rm CH}_{\bu} 
(\Xc_0)_{\Qb}
$$
(Proposition~1), with image $\wt{\rm ch} (\ov \Ec)$ in $\wt 
A (X)_{\Qb}$. Let also 
${\rm ch}_{\Yc_0} (Rf_* \, \Ec_{\Xc})
\in {\rm CH}_{\bu} (\Yc_0)_{\Qb}$ be the 
Chern character with supports in $\Yc_0$ of the triangle
$$
Rf_* \, S_{\Xc} \ra Rf_* \, E_{\Xc} \ra Rf_* \, Q_{\Xc} \ra 
Rf_* \, S_{\Xc} [1]
$$
(\S~2.4). 
Finally, let $\t_{\Yc'} (E_{\Xc})$, $\t_{\Yc'} (S_{\Xc})$ 
and $\t_{\Yc'} (Q_{\Xc})$ be defined as in (50) from 
$E_{\Xc}$, $S_{\Xc}$ and $Q_{\Xc}$ respectively. From (37) and 
(63) we must have
$$
\t_{\Yc'} (E_{\Xc}) - \t_{\Yc'} (S_{\Xc}) - \t_{\Yc'} 
(Q_{\Xc}) + p^* \, {\rm ch}_{\Yc_0} (Rf_* \, \Ec_{\Xc}) = 
f'_* \, \pi^* ({\rm ch}_{\Xc_0} (\Ec_{\Xc}) \, {\rm Td} 
(Tf)) \, . \leqno (68)
$$
Conversely, (67) and (68) will show that
 $f_* : {\build {K_0}_{}^{\vee}} (X) 
\ra {\build {K_0}_{}^{\vee}} (Y)$ exists, 
satisfying i), ii) 
and iii) in Theorem~1. After that, 
to check (65), we are  
reduced to the situation considered in ii). For any diagram 
(48) as in 4.1, the element ${\build {\rm ch}_{}^{\vee}} (\rho_* \, x)$ 
is then the projective system $(p^* \, {\rm ch} (Rf_* \, 
E_{\Xc}) + i_* \, \t_{\Yc'})$ in $\limind {\rm CH}^{\bu} 
(\Yc')_{\Qb}$. From Proposition~4, ii) we get
$$
p^* \, {\rm ch} (Rf_* \, E_{\Xc}) + i_* \, \t_{\Yc'} = f'_* 
\, \pi^* ({\rm ch} (E_{\Xc}) \, {\rm Td} (f)) \, ,
$$
which is precisely the $\Yc'$-component of
$$
\vp_* ({\build {\rm ch}_{}^{\vee}} (x) \, \wh{\rm Td} (\ov{T\vp})) \, .
$$
We are thus left with checking (67) and (68).

\medskip

\noindent {\bf 4.4.} Let us now
check the equality (67). By the 
definition ( ), if $K''$ is the cone of
$$
Ls^* \, Rf_* \, E_{\Xc} \ra Rf''_* \, L\s^* \, E_{\Xc} \, , 
\leqno (69)
$$
we have
$$
\leqalignno{
\t_{\Yc''} = \ &{\rm ch}_{\Yc''_0} (K) \, {\rm Td} (s) &(70) 
\cr
+ \ &s^* ({\rm ch} (Rf_* \, E_{\Xc})) \, {\rm Td}_0 (s) \cr
- \ &f''_* [\s^* ({\rm ch} (E_{\Xc}) \, {\rm Td} (f)) \, 
{\rm Td}_0 (\s)] \, , \cr
}
$$
and if $K'$ is the cone of
$$
Lr^* \, Rf'_* \, \pi^* \, E_{\Xc} \ra Rf''_* \, L\rho^* \, 
(\pi^* \, E_{\Xc})
$$
we have
$$
\leqalignno{
\t'_{\Yc''} = \ &{\rm ch}_{\Yc''_0} (K') \, {\rm Td} (r) &(71) 
\cr
+ \ &r^* ({\rm ch} (Rf'_* (\pi^* \, E_{\Xc}))) \, {\rm Td}_0 
(r) \cr
- \ &f''_* [\rho^* ({\rm ch} (\pi^* \, E_{\Xc}) \, {\rm Td} 
(f')) \, {\rm Td}_0 (\rho)] \, . \cr
}
$$
Since $E_{\Xc}$ is locally free we have
$$
L\rho^* (\pi^* \, E_{\Xc}) = \rho^* \, \pi^* \, E_{\Xc} = 
\s^* \, E_{\Xc} = L\s^* \, E_{\Xc} \, .
$$
Therefore the map (69) factors via $Lr^* \, Rf'_* (\pi^* \, 
E_{\Xc})$ and we get a triangle in $D_{\Yc''_0}^b (\Yc'')$
$$
Lr^* (K) \ra K'' \ra K \ra Lr^* (K) [1] \, .
$$
By (15) this implies
$$
{\rm ch}_{\Yc''_0} (K'') = r^* \, {\rm ch}_{\Yc'_0} (K) + 
{\rm ch}_{\Yc''_0} (K') \, . \leqno (72)
$$
Since ${\rm ch}_{\Yc''_0} (K'')$ is supported on $\Yc''_0$ 
we have (by (17))
$$
{\rm ch}_{\Yc''_0} (K'') \, {\rm Td} (s) = {\rm 
ch}_{\Yc''_0} (K'') \, {\rm Td}_0 (s) + {\rm ch}_{\Yc''_0} 
(K'') \, . \leqno (73)
$$
Similarly
$$
{\rm ch}_{\Yc''_0} (K') \, {\rm Td} (r) = {\rm ch}_{\Yc''_0} 
(K') \, {\rm Td}_0 (r) + {\rm ch}_{\Yc''_0} (K') \, . \leqno 
(74)
$$
Furthermore
$$
\leqalignno{
\ & s^* ({\rm ch} (Rf_* \, E_{\Xc})) \, {\rm Td}_0 (s) &(75) 
\cr
= \ & {\rm ch} (Ls^* \, Rf_* \, E_{\Xc}) \, {\rm Td}_0 (s) 
\cr
= \ & {\rm ch} (Rf''_* \, L\s^* \, E_{\Xc}) \, {\rm Td}_0 
(s) \cr
- \ & {\rm ch}_{\Yc''_0} (K'') \, {\rm Td}_0 (s) \cr
}
$$
and
$$
\leqalignno{
\ & r^* ({\rm ch} (Rf'_* (\pi^* \, E_{\Xc}))) \, {\rm Td}_0 
(r) &(76) \cr
= \ & {\rm ch} (Rf''_* \, L\s^* \, E_{\Xc}) \, {\rm Td}_0 
(r) \cr
- \ & {\rm ch}_{\Yc''_0} (K') \, {\rm Td}_0 (r) \, . \cr
}
$$

From (72)--(76) we conclude that
$$
\leqalignno{
\ & \t_{\Yc''} - \t'_{\Yc''} - r^* \, {\rm ch}_{\Yc''_0} (K) 
&(77) \cr
= \ & {\rm ch} (Rf''_* \, L\s^* \, E_{\Xc}) \, {\rm Td}_0 
(r) \cr
- \ & f''_* [\s^* ({\rm ch} (E_{\Xc}) \, {\rm Td} (f)) \, 
{\rm Td}_0 (\s)] \cr
- \ & {\rm ch} (Rf''_* \, L\s^* \, E_{\Xc}) \, {\rm Td}_0 
(r) \cr
+ \ & f''_* [\rho^* ({\rm ch} (\pi^* \, E_{\Xc}) \, {\rm Td} 
(f')) \, {\rm Td}_0 (\rho)] \, . \cr
}
$$
Applying the Riemann-Roch-Grothendieck formula to $f''$ and 
$\s^* \, E_{\Xc}$, we get from (77) that
$$
\t_{\Yc''} - \t'_{\Yc''} - r^* \, {\rm ch}_{\Yc''_0} (K) = 
f''_* ({\rm ch} (\s^* \, E_{\Xc}) \cdot A) \leqno (78)
$$
with
$$
\eqalign{
A := \ &{\rm Td} (f'') \, f''^* \, ({\rm Td}_0 (r) - {\rm 
Td}_0 (s)) \cr
- \ & \s^* ({\rm Td} (f)) \, {\rm Td}_0 (\s) \cr
+ \ & \rho^* ({\rm Td} (f')) \, {\rm Td}_0 (\rho) \, . \cr
}
$$
To compute $A$, we use the identities (44) and (45) from \S 3.1
to get
$$
\eqalign{
{\rm Td} (f'') \, f''^* \, {\rm Td}_0 (r) = & \ \wt{\rm Td} 
(rf'' , f'') \, , \cr
{\rm Td} (f'') \, f''^* \, {\rm Td}_0 (s) = & \ \wt{\rm Td} 
(sf'' , f'') \, , \cr
\s^* ({\rm Td} (f)) \, {\rm Td}_0 (\s) = & \ \wt{\rm Td} 
(f\s ,f) \, , \cr
\rho^* ({\rm Td} (f')) \, {\rm Td}_0 (\rho) = & \ \wt{\rm 
Td} (f'\rho ,f') \, . \cr
}
$$
Since $sf'' = f\s$ and $rf'' = f'\rho$, we get from this and 
(28)
$$
\leqalignno{
A = \ & \wt{\rm Td} (sf'' , f'') - \wt{\rm Td} (rf'' , f'') 
&(79) \cr
- \ & \wt{\rm Td} (f\s , f) + \wt{\rm Td} (f' \rho , f') \cr
= \ & \wt{\rm Td} (f'' , f) - \wt{\rm Td} (f'' , f') = 
\rho^* \, \wt{\rm Td} (f',f)_{\Yc'_0} \, . \cr
}
$$
From (78) and (79) the equality (67) follows.

\medskip

\noindent {\bf 4.5.} To check (68) in \S~4.3 we let 
$K(E_{\Xc})$ be the cone of
$$
Lp^* \, Rf_* \, E_{\Xc} \ra Rf'_* \, L\pi^* \, E_{\Xc}
$$
and we define  $K(S_{\Xc})$ and $K(Q_{\Xc})$ similarly. Since 
these maps fit in a morphism of triangles
$$
\diagram{
Lp^* \, Rf_* \, S_{\Xc} &\hfl{}{} &Lp^* \, Rf_* \, E_{\Xc} 
&\hfl{}{} &Lp^* \, Rf_* \, Q_{\Xc} &\hfl{}{} &Lp^* \, Rf_* 
\, S_{\Xc} [1] \cr
\vfl{}{} &&\vfl{}{} &&\vfl{}{} &&\vfl{}{} \cr
Rf'_* \, L\pi^* \, S_{\Xc} &\hfl{}{} &Rf'_* \, L\pi^* \, 
E_{\Xc} &\hfl{}{} &Rf'_* \, L\pi^* \, Q_{\Xc} &\hfl{}{} 
&Rf'_* \, L\pi^* \, S_{\Xc} [1] \, , \cr
}
$$
where all lines and columns are acyclic on $Y$, we get from 
(30) in Lemma~1 i) that
$$
\leqalignno{
\ & {\rm ch}_{\Yc'_0} (K (E_{\Xc})) - {\rm ch}_{\Yc'_0} (K 
(S_{\Xc})) - {\rm ch}_{\Yc'_0} (K (Q_{\Xc})) &(80) \cr
= \ & {\rm ch}_{\Yc'_0} (Rf'_* \, L\pi^* \, \Ec_{\Xc}) - 
{\rm ch}_{\Yc'_0} (Lp^* \, Rf_* \, \Ec_{\Xc}) \, , \cr
}
$$
where $Rf'_* \, L\pi^* \, \Ec_{\Xc}$ is the upper triangle 
and
$$
{\rm ch}_{\Yc'_0} (Lp^* \, Rf_* \, \Ec_{\Xc}) = p^* \, {\rm 
ch}_{\Yc_0} (Rf_* \, \Ec_{\Xc}) \, . \leqno (81)
$$
On the other hand,
$$
\leqalignno{
\ & p^* \, [{\rm ch} (Rf_* \, E_{\Xc}) - {\rm ch} (Rf_* \, 
S_{\Xc}) - {\rm ch} (Rf_* \, Q_{\Xc})] \, {\rm Td}_0 (p) &(82) 
\cr
= \ & p^* \, {\rm ch}_{\Yc_0} (Rf_* \, \Ec_{\Xc}) \, {\rm 
Td}_0 (p) \cr
= \ & p^* \, {\rm ch}_{\Yc_0} (Rf_* \, \Ec_{\Xc}) \, {\rm 
Td} (p) - p^* \, {\rm ch}_{\Yc_0} (Rf_* \, \Ec_{\Xc}) \, . 
\cr
}
$$
Since $f \pi = p f'$ we have
$$
\pi^* ({\rm Td} (f)) \, {\rm Td} (\pi) = {\rm Td} (f') \, 
f'^* ({\rm Td} (p)) \, ,
$$
and this implies that
$$
\leqalignno{
\ & f'_* [\pi^* (({\rm ch} (E_{\Xc}) - {\rm ch} (S_{\Xc}) - 
{\rm ch} (Q_{\Xc})) \, {\rm Td} (f)) \, {\rm Td}_0 (\pi)] 
&(83) \cr
= \ & f'_* [\pi^* ({\rm ch}_{\Xc_0} (\Ec_{\Xc}) \, {\rm Td} 
(f)) \, {\rm Td}_0 (\pi)] \cr
= \ & f'_* [\pi^* ({\rm ch}_{\Xc_0} (\Ec_{\Xc}) \, \pi^* \, 
{\rm Td} (f) \, {\rm Td} (\pi)] \cr
\ & - f'_* \, \pi^* ({\rm ch}_{\Xc_0} (\Ec_{\Xc}) \, {\rm 
Td} (f)) \cr
= \ & f'_* [\pi^* ({\rm ch}_{\Xc_0} (\Ec_{\Xc})) \, {\rm Td} 
(f') \, f'^* \, {\rm Td} (p)] \cr
\ & - f'_* \, \pi^* ({\rm ch}_{\Xc_0} (\Ec_{\Xc}) \, {\rm 
Td} (f)) \cr
= \ & {\rm ch}_{\Yc'_0} (Rf'_* \, L\pi^* \, \Ec_{\Xc}) \, 
{\rm Td} (p) \cr
\ & - f'_* \, \pi^* ({\rm ch}_{\Xc_0} (\Ec_{\Xc}) \, {\rm 
Td} (f)) \, , \cr
}
$$
where the last equality follows from the projection formula 
for $f'$ together with the Riemann-Roch-Grothen\-dieck 
formula with supports. Putting (50), (80), (81), (82) and (83) 
together, we get
$$
\t_{\Yc'} (E_{\Xc}) - \t_{\Yc'} (S_{\Xc}) - \t_{\Yc'} 
(Q_{\Xc}) + p^* \, {\rm ch}_{\Yc_0} (Rf_* \, \Ec_{\Xc}) = 
f'_* \, \pi^* ({\rm ch}_{\Xc_0} (\Ec_{\Xc}) \, {\rm Td} (f)) 
\, ,
$$
i.e. (68) holds true. This ends the proof of Theorem~1.

\medskip

\noindent {\bf 4.6.} Here are a few more properties of the 
direct image morphisms defined in Theorem~1.

\medskip

\noindent {\bf Theorem 2.} {\it Let $f : X \ra Y$ be a map 
of varieties, equipped with an arbitrary metric on $Tf$.}

\item{i)} {\it When $x \in {\build {K_0}_{}^{\vee}} (X)$ and $y \in 
{\build {K_0}_{}^{\vee}} (Y)$, we have }
$$
f_* (xf^* (y)) = f_* (x) y \, .
$$

\item{ii)} {\it If $f$ is flat, $f_*$ maps $\wh{K}_0 (X)$ 
into $\wh{K}_0 (Y)$.}

\item{iii)} {\it Let $g : Y \ra Z$ be a map of varieties. 
Choose arbitrary metrics on $Tg$ and $T(gf)$. Let $\wt{\rm 
Td}$ be the secondary Todd class of the metrized triangle on 
$X$
$$
\Tc : Tf \ra T(gf) \ra f^* Tg \ra Tf [1] \, .
$$
Then, for any $x \in {\build {K_0}_{}^{\vee}} (X)$, the following 
identity holds:}
$$
(gf)_* (x) - g_* (f_* (x)) = - \a ((gf)_* ({\rm ch} (x) \, 
\wt{\rm Td})) \, . \leqno (84)
$$

\smallskip

\noindent {\bf 4.7.} To check Theorem~2 i) when $x$ or $y$ 
is in the image of $\a$, we just use the fact that the 
projection formula is true for forms and currents (\S~1.8 
and [BGS] 1.5), together with (63).

\smallskip

We may then assume that $x$ is the class of $(E,h,0)$ and 
$y$ is the class of $(F,h',0)$. By the previous argument, 
and the anomaly formula (66), the difference
$$
z = f_* (x \, f^* (y)) - f_* (x) \, y
$$
does not depend on the choice of the metrics $h$ and $h'$. 
Furthermore, if we change the metric on $Tf$ and if we
denote by 
$\wt{\rm Td}$ the corresponding secondary Todd class, it 
follows from Theorem~1,  iii) that $z$ gets 
replaced by
$$
z + \a \, f_* ({\rm ch} (x \, f^* (y)) \, \wt{\rm Td}) - \a 
(f^* ({\rm ch} (x) \, \wt{\rm Td})) \cdot y \, .
$$
This is equal to $z$ since, by (39),
$$
\a (f_* ({\rm ch} (x) \, \wt{\rm Td})) \cdot y = \a (f_* ({\rm 
ch} (x) \, \wt{\rm Td}) \, {\rm ch} (y)) = \a \, f_* ({\rm 
ch} (x \, f^* (y)) \, \wt{\rm Td}) \, .
$$
Consequently we may assume that there is a map of models 
$\ud f : \Xc \ra \Yc$ inducing $f$ on $X$ and defining the 
metric on $Tf$, and that the metric on $E$ (resp. $F$) is defined 
by a bundle $E_{\Xc}$ (resp. $F_{\Yc}$) on $\Xc$ (resp. 
$\Yc$). Since $Rf_* (E_{\Xc} \ot f^* \, F_{\Yc})$ is 
isomorphic to $Rf_* \, E_{\Xc} \ot F_{\Yc}$ we know that $z$ 
can be written
$$
z = \a (\eta) \, , \ \eta \in \wt D (Y)_{\Qb} \, .
$$
On the other hand, the Riemann-Roch formula (65) in Theorem~1, 
 together with the projection formula for arithmetic Chow 
groups (\S~1.8) imply that $\wh{\rm ch} (z) = 0$. Therefore
$$
dd^c \, \eta = \om (\wh{\rm ch} (z)) = 0 \, .
$$
By Proposition~1, ii), if $\eta_{\Yc}$ vanishes we can 
conclude that $\eta = 0$ and $z=0$. But the analytic torsion 
$\t_{\Yc}$ for both $E_{\Xc}$ and $E_{\Xc} \ot f^* \, 
F_{\Yc}$ are zero by (51). Therefore $\eta_{\Yc} = 0$.

\medskip

\noindent {\bf 4.8.} Assume $f$ is flat and $x \in \wh{K}_0 
(Y)$. Since ${\rm ch} (x) \, {\rm Td} (\ov{Tf})$ lies in 
$A_{\rm closed} (X)_{\Qb}$ and $f_*$ maps forms into forms 
(\S~1.7, i.e. [BGS] (4.1.1) and (4.2.1)) $f_* ({\rm ch} (x) 
\, {\rm Td} (\ov{Tf}))$ lies in $A_{\rm closed} (Y)_{\Qb}$. 
But it follows from 
the Riemann-Roch formula (65) that
$$
f_* ({\rm ch} (x) \, {\rm Td} (\ov{Tf})) = {\rm ch} (f_* (x)) \, .
$$
Therefore $f_* (x)$ lies in $\wh{K}_0 (Y)$ (see 2.5).

\medskip

\noindent {\bf 4.9.} The proof of Theorem~2, iii) is similar 
to Th.~2, i). Namely, let
$$
z = (gf)_* \, (x) - g_* (f_* (x)) + \a ((gf)_* \, ({\rm ch} 
(x) \, \wt{\rm Td})) \, .
$$
When $x = \a (\eta)$, it follows from (63) that 
$z=0$. Indeed, since, by (25) and (43),
 $${\rm Td} (\ov{T(gf)}) = {\rm Td} 
(\ov{Tf}) \, f^* \, {\rm Td} (\ov{Tg}) - dd^c (\wt{\rm 
Td}) \, ,$$ we get
$$
\eqalign{
(gf)_* \, (x) = \ &\a (g_* (f_* (\eta \, {\rm Td} (\ov{Tf})) 
\, {\rm Td} (\ov{Tg})) \cr
\ & - \a ((gf)_* \, (\eta \, dd^c (\wt{\rm Td}))) \cr
= \ & g_* \, f_* (x) - (gf)_* \, (\a (\eta \, dd^c (\wt{\rm 
Td}))) \cr
}
$$
and $z=0$ by the ``Stokes formula''
$$
\eta \, dd^c (\om) = dd^c (\eta) \, \om \, .
$$

Therefore it is enough to check that $z=0$ when $x$ is the 
class of $(E,E_{\Xc},0)$ for some extension of $E$ to a 
model of $X$. When the metric $h_f$ on $Tf$ is replaced by 
$h_{f'}$, from Theorem~1, iii) and (63) we know that $z$ gets 
replaced by
$$
z' = z - \a g_* [f_* ({\rm ch} (x) \, \wt{\rm Td} (h_f , 
h'_f)) \, {\rm Td} (\ov{Tg})] - \a (gf)_* [{\rm ch} (x) 
(\wt{\rm Td} - \wt{\rm Td}')] \, ,
$$
where $\wt{\rm Td}$ (resp. $\wt{\rm Td}'$) are the secondary 
Todd classes of the triangle $\Tc$, where the metric on $Tf$ 
is $h_f$ (resp. $h'_f$), and we do not change the metrics on 
$T(gf)$ and $f^* \, Tg$. From Lemma~1 ii) we know that
$$
\wt{\rm Td} - \wt{\rm Td}' = \wt{\rm Td} (h_f , h'_f) \, f^* 
\, {\rm Td} (\ov{Tg}) \, ,
$$
therefore $z' = z$.

\smallskip

One checks in a similar way (by shifting (20)) that $z$ does 
not depend on the metrics on $Tg$ and $T(gf)$. Finally one 
may assume that there are maps of models $\ud f : \Xc \ra  
\Yc$ and $\ud g : \Yc \ra \Zc$ which induce $f$ and $g$, that 
the metric on $E$ is given by a bundle $E_{\Xc}$ on $\Xc$, 
and that the metrics on $Tf$ (resp. $Tg$, resp. $T(gf)$) is 
given by $T \ud f$ (resp. $T \ud g$, resp. $T (\ud g \ud 
f)$). The exact triangle 
$$T \ud f \ra T (\ud g \ud f) \ra 
\ud f^* \, T \ud g \ra T \ud f [1]$$
 implies in that case 
that $\wt{\rm Td} = 0$. Since $R \ud f_* \, R \ud g_* = R 
(\ud f \ud g)_*$ we can write $z = \a (\eta)$, $\eta \in \wt 
D (Z)_{\Qb}$, and from Riemann-Roch we deduce that
$$
dd^c (\eta) = {\rm ch} (z) = (gf)_* \, ({\rm ch} (x) \, {\rm 
Td} (\ov{T(gf)})) - g_* (f_* ({\rm ch} (x) \, {\rm Td} 
(\ov{Tf})) \, {\rm Td} (\ov{Tg})) = 0 \, .
$$
To check that $\eta = 0$ (hence $z=0$) it is enough that 
$\eta_{\Zc} = 0$. But this follows from (51) applied to the maps 
$\ud f$, $\ud g$ and $\ud g \ud f$. This ends the proof of 
Theorem~2.

\bigskip
\bigskip
\bigskip

\centerline{\bf References}

\medskip

 \item{[AKMW]}
D. Abramovich, K. Karu, K. Matsuki, J. Wlodarczyk :
 Torification and factorization of birational maps,
preprint, 1999, {\bf math.AG/9904135}.

\smallskip

 \item{[BFM]}
P. Baum, W. Fulton, R. MacPherson
: Riemann-Roch for Singular Varieties
{\it Pub. Math. I.H.E.S.}
{\bf 45}, 1975,
 253--290.

\smallskip

 \item{[BK]}
J.-M. Bismut, K. Koehler:
Higher analytic torsion forms 
for direct images and anomaly formulas,
{\it J. Algebr. Geom.} {\bf 1} No.4, 1992, 647-684.

\smallskip

\item{[BGS]} S. Bloch, H. Gillet, C. Soul\'e : 
Non-archimedean Arakelov theory, {\it Journal of Algebraic
Geometry}, {\bf 4}, 1995,  427-485.

\smallskip

 \item{[B]} J. Burgos:
Arithmetic Chow rings and Deligne-Beilinson cohomology,
{\it J. Algebr. Geom.} {\bf 6} No.2, 1997, 335-377.
 
\smallskip

\item{[F]} W. Fulton : Intersection theory, {\it Ergebnisse
der Math.} {\bf 3}, Folge 2 Band 2, 1984, Springer-Verlag,
Berlin-Heidelberg-New York.

\smallskip

\item{[Fr]} J. Franke : Riemann-Roch in functorial form,
preprint, 1992, 78 pp.

\smallskip

 \item{[GS1]} H.Gillet, C.Soul\'e : Arithmetic Intersection Theory,
{\it Publications Math. IHES} {\bf 72}, 1990, 94-174.

\smallskip

 \item{[GS2]} H.Gillet, C.Soul\'e :
 Characteristic classes for algebraic vector bundles with
hermitian metric, {\it Annals of Maths.} {\bf  131}, 1990,
163-203.

\smallskip

 \item{[GS3]} H.Gillet, C.Soul\'e :
 Analytic torsion and the Arithmetic Todd genus,
{\it Topology}
{\bf 30}, 1 ,1991, 21-54.

\smallskip

 \item{[GS4]} H.Gillet, C.Soul\'e :
 An arithmetic Riemann-Roch theorem,
{\it Inventiones Math.} {\bf  110}, 1992, 474-543.

 \item{[H]} H. Hironaka: Resolution of singularities of an algebraic
variety over a field of characteristic zero
, {\it Annals of Math.}
{\bf 79},
 1964
, 109--326.

\smallskip

 \item{[M]} H. Matsumura: 
Commutative ring theory. 
Transl. from the Japanese by M. Reid,
 {\it Cambridge Studies in Advanced Mathematics},
 {\bf 8},  Cambridge University Press,
 1989. 

\smallskip

 \item{[Q]}
  D. Quillen:  Determinants of Cauchy--Riemann operators over a
Riemann surface, {\it Funct. Anal. Appl.\/}, 1985, 31-34.

\smallskip

 \item{[RG]} M. Raynaud,  L. Gruson:
Crit\`{e}res de platitude et de projectivit\'{e}
{\it Inv. Math.}
{\bf  13}
, 1971,
 1--89.
 
 \smallskip

 \item{[S]} T. Saito:
Conductor, discriminant, and the Noether formula 
of arithmetic surfaces
{\it Duke Math. Journal}
{\bf  57}
, 1988,
 151--173.
\smallskip

 \item{[SGA4]}
  M. Artin, A. Grothendieck, J.L. Verdier, P. Deligne, B. Saint-Donat:
S\'eminaire de g\'eom\'etrie alg\'ebrique du Bois-Marie 1963-1964, 
Th\'eorie des topos et cohomologie \'etale des sch\'emas ,SGA 4,
 Tome 3, Expos\'es IX a XIX. 
{\it Lecture Notes in Mathematics} {\bf 305}, 
1973, Berlin-Heidelberg-New York: Springer-Verlag.

\smallskip

 \item{[SGA6]}
P. Berthelot, A. Grothendieck, L. Illusie :
S\'eminaire de g\'eom\'etrie alg\'ebrique du Bois Marie 
1966/67, SGA 6, 
Th\'eorie des intersections et th\'eor\`eme de Riemann-Roch 
,{\it Lecture Notes in Mathematics } {\bf 225}, 1971,
 Berlin-Heidelberg-New York: Springer-Verlag. 
 
 \smallskip

 \item{[W]}
 J. Wlodarczyk : Combinatorial structures on toroidal varieties
and a proof of the weak factorization theorem
preprint, 1999, {\bf math.AG/9904076}.

 \smallskip

 \item{[Z]} J. Zha : A general Arithmetic Riemann-Roch theorem, PHD
 thesis, Chicago University, 1998.

\bigskip

H.G.:  Department of Mathematics, Statistics, and Computer Science,
 University of Illinois at Chicago, 851 S. Morgan Street,
Chicago, IL 60607-7045, U.S.A.

\medskip

C.S. : CNRS, Institut des Hautes \'{E}tudes Scientifiques,  35, Route de
Chartres, 91440, Bures-sur-Yvette, France

\bye